\documentclass[times]{article}

\usepackage{moreverb}
\usepackage[margin=1.25in]{geometry}
\usepackage[dvips,colorlinks,bookmarksopen,bookmarksnumbered,citecolor=red,urlcolor=red]{hyperref}
\usepackage{graphicx}
\usepackage{amssymb}
\usepackage{amsmath,bm}
\usepackage{color}
\usepackage{xspace}
\usepackage{algpseudocode}
\usepackage{algorithmicx}
\usepackage{subfig}
\usepackage{url}
\usepackage{comment}
\usepackage[subnum]{cases}
\usepackage{psfrag}

\usepackage{tikz}
\usetikzlibrary{shapes,arrows}
\usetikzlibrary{positioning}
\tikzstyle{block}=[draw opacity=0.7,line width=1.4cm]
\definecolor{CranJ}{cmyk}{0,0.69,0.54,0.04} 
\definecolor{PinkJ}{cmyk}{0,0.71,0.43,0.12} 
\definecolor{Cran}{cmyk}{0,0.73,0.41,0.29} 
\definecolor{VRed}{cmyk}{0,0.75,0.25,0.2} 
\definecolor{ORed}{cmyk}{0,0.75,0.75,0} 
\definecolor{CBlue}{cmyk}{1,0.25,0,0} 

\newcommand\BibTeX{{\rmfamily B\kern-.05em \textsc{i\kern-.025em b}\kern-.08em
T\kern-.1667em\lower.7ex\hbox{E}\kern-.125emX}}

\newcommand{\VV}{\mathcal{V}}
\newcommand{\EE}{\mathcal{E}}
\newcommand{\GG}{\mathcal{G}}
\newcommand{\PP}{\mathcal{P}}

\newcommand{\real}{{\mathbb{R}}}

\newcommand{\Hlambda}{\hat{\lambda}}

\newcommand{\until}[1]{\{1,\dots,#1\}}
\newcommand{\Lnorm}{\Big\|}
\newcommand{\Rnorm}{\Big\|}

\newcommand{\overbar}[1]{\mkern 1.5mu\overline{\mkern-1.5mu#1\mkern-1.5mu}\mkern 1.5mu}
\newcommand{\barshift}[1]{\skew{5}{\overbar}{#1}}
\newcommand{\tildeshift}[1]{\skew{5}{\widetilde}{#1}}

\newcommand{\vect}[1]{\boldsymbol{#1}}
\newcommand{\Bvect}[1]{\barshift{\bm{{#1}}}}
\newcommand{\Tvect}[1]{\tildeshift{\bm{{#1}}}}
\newcommand{\vectT}[1]{\vect{#1}_{\mathrm{T}}}
\newcommand{\dvect}[1]{\dot{\vect{#1}}}
\newcommand{\ddvect}[1]{\ddot{\vect{#1}}}

\newcommand{\e}[1]{\operatorname{e}^{#1}}
\newcommand{\Sym}[1]{\operatorname{Sym}(#1)}
\newcommand{\Diag}[1]{\operatorname{Diag}(#1)}
\newcommand{\sat}[2]{\operatorname{sat}_{#1}(#2)}

\newcommand{\avrg}[1]{\frac{1}{N}\sum_{j=1}^N#1}
\newcommand{\SUM}[2]{\sum_{#1}^{#2}}
 \newcommand{\boxend}{\hfill \ensuremath{\Box}}

\definecolor{PineGreen} {cmyk}{0.92,0,0.59,0.25}

\newtheorem{thm}{Theorem}[section]
\newtheorem{pro}{Proposition}[section]
\newtheorem{rem}{Remark}[section]
\newtheorem{cor}{Corollary}[section]
\newtheorem{lem}{Lemma}[section]
\newtheorem{defn}{Definition}

\newtheorem{prob}{Problem}
\makeatletter
\renewcommand*{\@opargbegintheorem}[3]{\trivlist
      \item[\hskip \labelsep{\emph{ #1\ #2}}] \emph{(#3):}\ \itshape}
\makeatother

\begin{document}

\date{}

\title{Dynamic Average Consensus
 under Limited Control Authority and Privacy Requirements\footnote{under review in International Journal of Robust and Nonlinear Control}}

\author{Solmaz~S.~Kia, Jorge Cort\'es, Sonia Mart{\'\i}nez\\
{\large Department of Mechanical and Aerospace Engineering,}\\
 {\large University of California San Diego, La Jolla, CA 92093, USA}}

\maketitle
\begin{abstract}
  This paper introduces a novel continuous-time dynamic average
  consensus algorithm for networks whose interaction is described by a
  strongly connected and weight-balanced directed graph.  The proposed
  distributed algorithm allows agents to track the average of their
  dynamic inputs with some steady-state error whose size can be
  controlled using a design parameter.  This steady-state error
  vanishes for special classes of input signals.  We analyze the
  asymptotic correctness of the algorithm under time-varying
  interaction topologies and characterize the requirements on the
  stepsize for discrete-time implementations.  We show that our
  algorithm naturally preserves the privacy of the local input of each
  agent.  Building on this analysis, we synthesize an extension of the
  algorithm that allows individual agents to control their own rate of
  convergence towards agreement and handle saturation bounds on the
  driving command.  Finally, we show that the proposed extension
  additionally preserves the privacy of the transient response of the
  agreement states and the final agreement value from internal and
  external adversaries. Numerical examples illustrate the results.
\end{abstract}
 
\section{Introduction}
\vspace{-2pt}

This paper studies the dynamic average consensus problem for a network
of autonomous agents. Given a set of time-varying signals, one per
agent, this problem consists of designing a distributed algorithm that
allow agents to track the time-varying average of the signals using
only information from neighbors.  Solutions to this problem are of
interest in scenarios that require the fusion of dynamic and evolving
information collected by multiple agents. Examples include multi-robot
coordination~\cite{PY-RAF-KML:08}, distributed spatial
estimation~\cite{SM:07e,JC:07-tac}, sensor
fusion~\cite{ROS-JSS:05,ROS:07}, feature-based map
merging~\cite{RA-JC-CS:10j}, and distributed
tracking~\cite{PY-RAF-KML:07}. We are particularly interested in
algorithmic solutions that allow agents to adjust the rate of
convergence towards agreement, are able to handle constraints on
actuation, and preserve the privacy of the information available to
them against adversaries. 

\vspace{-6pt}
\subsubsection*{Literature review.}
\vspace{-2pt} 

Consensus problems have been intensively studied over the last years.
The main body of work focuses on the static case, where agents aim to
reach consensus on a function depending on initial static values, see
e.g.~\cite{ROS-JAF-RMM:07,WR-RWB:08,WR-YC:11,MM-ME:10,FB-JC-SM:09} and
references therein. In contrast, the literature on dynamic consensus
is not as rich.  The initial work~\cite{DPS-ROS-RMM:05b} proposes a
dynamic average consensus algorithm that under proper initialization
is able to track, with zero steady-state error, the average of dynamic
inputs whose Laplace transfer functions have at most one pole at the
origin and the rest of the poles are in the left half-plane.
In~\cite{ROS-JSS:05}, the
authors generalize the static consensus algorithm of~\cite{ROS-RMM:04}
to track the average of inputs with bounded
derivatives 
which differ by a zero-mean Gaussian noise. The algorithm acts as a
low-pass filter that allows agents to track the average of dynamic
inputs with a non-zero steady-state error, which vanishes in the
absence of noise. Using input-to-state stability analysis,~\cite{RAF-PY-KML:06} proposes
a proportional-integral algorithm to solve the dynamic consensus
problem which, from any initial condition, converges with non-zero
steady-state error if the signals are slowly time-varying, and exactly
if the signals are static. This algorithm is generalized
in~\cite{HB-RAF-KML:10} to achieve zero-error dynamic average
consensus of a special class of time-varying input signals whose
Laplace transform is a rational function with no poles in the
left-hand complex plane.  The proposed algorithm employs
frequency-domain tools and exploits the properties of the inputs'
Laplace transforms.  All the algorithms mentioned above are designed
in continuous time and work for networks with a fixed, connected, and
undirected graph topology. The results of~\cite{RAF-PY-KML:06} can be
applied to networks with a strongly connected and weight-balanced
digraph topology provided each agent can communicate with its
out-neighbors and knows the weights of its incoming edges. Such
requirement may be hard to satisfy in scenarios where the topology is
changing.  The work~\cite{MZ-SM:08a} develops an alternative class of
discrete-time dynamic average consensus algorithms  whose
convergence analysis relies on input-to-output stability properties in
the presence of external disturbances.  With a proper initialization
of the states, the proposed schemes can track, with a bounded
steady-state error, the average of the time-varying inputs whose
$n$th-order difference is bounded. If the $n$th-order difference is
asymptotically zero, the estimates of the average converge to the true
average asymptotically with one timestep delay.    
Other classes of algorithms related to our work are leader-follower
algorithms for networks of mobile agents with integrator dynamics,
e.g., see~\cite{WR:07a,GS-YH-KHJ:12}, and robust average consensus
algorithms in the presence of additive input
disturbances~\cite{GS-KHJ:13}. In the former scenario, agents reach
consensus by following the input signal of the leader agent(s),
instead of converging to the average of input signals across the
network. In the latter case, the algorithm performance achieving
consensus is analyzed in the presence of dynamic external
disturbances.
A common limitation of the works cited above is the lack of
consideration of restrictions on the rate of convergence of individual
agents, bounded control authority, or privacy issues. Regarding the
latter, the above algorithms require agents to share their agreement
state with their neighbors, and, in some cases, even their local
inputs. Therefore, if adversaries are able to listen to the exchanged
messages, they could infer local inputs, sensitive transient responses
and final agreement states of the network.

\vspace{-6pt}
\subsubsection*{Statement of contributions.}
\vspace{-2pt}%

We begin by providing a formal statement of the dynamic average
consensus problem for a multi-agent system, paying special attention
to the rate of convergence, limits on control actuation, and the
preservation of privacy.  Our starting point is the introduction of a
continuous-time algorithm that allows the group of agents
communicating over a strongly connected and weight-balanced digraph to
track the average of their reference inputs with some steady-state
error.  We carefully characterize the asymptotic convergence
properties of the proposed strategy, including its rate of
convergence, its robustness against initialization errors, and its
amenability to discrete-time implementations. We also discuss how the
algorithm performance (specifically, the steady-state error and the
transient response) can be tuned via two design parameters.  For
special classes of inputs, which include static inputs and dynamic
inputs which differ by a constant value, we show that the steady-state
error vanishes. We also establish the algorithm correctness under
time-varying network topologies that remain weight-balanced and are
infinitely often jointly strongly connected.  Our next step is the
introduction of an extension of the proposed dynamic average consensus
algorithm to include a local first-order filter at each agent. We show
how this extension allows individual agents to tune their rate of
convergence towards agreement without affecting the rest of the
network or changing the ultimate tracking error bound. We also
establish that, under limited control authority, this extension has
the same correctness guarantees as the original algorithm as long as
the input signals are bounded with a bounded relative growth. Several
simulations illustrate our results.
Our final step is the characterization of the privacy-preservation
properties of the proposed dynamic average consensus algorithms. We
consider adversaries who aim to retrieve information about the inputs,
their average, or the state trajectories. These adversaries might be
inside (internal) or outside (external) the network, do not interfere
with the algorithm execution, and may have access to different levels
of information, such as knowledge of certain parts of the graph
topology, the algorithm design parameters, initial conditions, or the
history of communication messages.  We show how the proposed
algorithms naturally preserve the privacy of the input of each agent
against any adversary.  Moreover, we establish that the extension that
incorporates local first-order filters protects the privacy of the
agreement state trajectories against any adversary by adding a common
signal to the messages transmitted among neighbors. This strategy also
preserves the privacy of the final agreement value against external
adversaries.

\vspace{-6pt}
\subsubsection*{Organization.}
\vspace{-2pt} Section~\ref{sec::prelim} introduces basic notation,
graph-theoretic concepts, and the model of time-varying
networks. Section~\ref{sec::ProbDef} formally introduces the dynamic
consensus problems of interest.  Section~\ref{sec::DyConsensus}
presents our dynamic average consensus algorithm, establishes its
correctness, and analyzes its properties regarding changing
interaction topologies, discrete-time implementations, and rate of
convergence. Section~\ref{sec::Rate} introduces a modified version
which enables agents to opt for a slower rate of convergence and
solves the consensus problem in the presence of bounded control
commands.  Section~\ref{sec::Priv} considers the privacy preservation
properties of the proposed algorithms.  Section~\ref{sec:simulations}
presents simulations illustrating our results.  Finally,
Section~\ref{sec:conclusions} gathers our conclusions and ideas for
future work.

\vspace{-6pt}
\section{Preliminaries}\label{sec::prelim}\vspace{-2pt}
In this section, we introduce basic notation, concepts from graph
theory used throughout the paper, and our model for networks with
time-varying interaction topologies.

\vspace{-6pt}
\subsection{Notational conventions}\vspace{-2pt}
The vector $\vect{1}_n$ is the vector of $n$ ones, $\vect{0}_n$ is the
vector of $n$ zeros, and $\vect{I}_n$ is the identity matrix with
dimension $n\times n$. We denote by $\vect{A}^\top$ the transpose of
matrix~$\vect{A}$. For a square matrix $\vect{A}$ we define
$\Sym{\vect{A}}\ = \frac{1}{2} (\vect{A} +\vect{A}^\top )$. We use
$\mathrm{Diag}(\vect{A}_1,\cdots,\vect{A}_N)$ to represent the
block-diagonal matrix constructed from matrices
$\vect{A}_1,\dots,\vect{A}_N$.  We define
$\vect{\Pi}_n=\vect{I}_n-\frac{1}{n}\vect{1}_n\vect{1}_n^\top$.
We denote the induced two-norm of a real matrix $\vect{A}$ by
$\|\vect{A}\|$, i.e., $\|\vect{A}\|=\sigma_{\max}(\vect{A})$, where
$\sigma_{\max}$ is the maximum singular value of~$\vect{A}$.  The
spectral radius of a square matrix $\vect{A}$ is represented by
$\rho(\vect{A})$. For a vector $\vect{u}$, we use $\|\vect{u}\|$ to
denote the standard Euclidean norm, i.e.,
$\|\vect{u}\|=\sqrt{\vect{u}^\top\vect{u}}$.  For vectors
$\vect{u}_1,\cdots,\vect{u}_N$, we let
$(\vect{u}_1,\cdots,\vect{u}_N)$ represent their aggregated
vector. For a complex variable $c$, $\Re(c)$ indicates its real part.
For a scalar variable $u$, the saturation function with limit
$0<\bar{u}<\infty$ is indicated by $\sat{\bar{u}}{u}$,
i.e., $\sat{\bar{u}}{u}=\text{sign}(u)\min\{|u|, \bar{u}\}$.  We let
$\delta_1(\epsilon)\in O(\delta_2(\epsilon))$ denote the fact that
there exist positive constants $c$ and $k$ such that $
|\delta_1(\epsilon)|\leq k|\delta_2(\epsilon)|,~\forall
\,|\epsilon|<c$. For network-related variables, the local variables of
each agent are distinguished by a superscript, e.g., $u^i(t)$ is the
local dynamic input of agent $i$. If $p^i\in\real$ is a local variable
at agent $i$, the aggregated $p^i$'s are represented by $\vect{p} =
(p^1,\dots, p^N) \in \real^{N}$.  Our analysis involves linear systems
of the form
\begin{equation}\label{eq::GModel}
  \dvect{x}(t)=\vect{A}\vect{x}(t)+\vect{B}\vect{u}(t),
\end{equation}
where states $\vect{x}(t)$ take values in the Euclidean space
$\real^n$, and inputs are measurable locally essentially bounded maps
$\vect{u} : [0,\infty)\to\real^{m}$. The \emph{zero-system} associated
to \eqref{eq::GModel} is by definition the system with no inputs,
i.e., $\dvect{x}=\vect{A}\vect{x}$.  We denote by
$\|\vect{u}\|_{\text{ess}}$, the (essential) supremum norm , i.e.,
$\|\vect{u}\|_{\text{ess}} = \sup\{\|\vect{u}(t)\|,~t\geq0\}<\infty$.
The \emph{convergence rate} of a stable linear system
$\dvect{x}=\vect{A}\vect{x}$ is
\begin{align}\label{eq::rate_of_conv}
  r &= \inf\{\chi>0\,| \,\exists\, \kappa>0 \text{ such that }
  \|\vect{x}(t)\|\leq \kappa\|\vect{x}(0)\|\e{-\chi t},~~t\geq 0\}.
\end{align}
Here, $\vect{x}(t)$ is the solution of the system when it starts from
any initial state $\vect{x}(0)\in\real^n$.  This definition implies
that for a linear time-invariant dynamical system, the rate of
convergence is the least negative real part of the eigenvalues of the
system matrix.

\vspace{-6pt}
\subsection{Graph theory}\label{app::Graph}\vspace{-2pt}

We briefly review some basic concepts from the graph, see
e.g.~\cite{FB-JC-SM:09}.  A \emph{directed graph}, or simply a
\emph{digraph}, is a pair $\GG = (\VV ,\EE )$, where
$\VV=\{1,\dots,N\}$ is the \emph{node set} and $\EE \subseteq
\VV\times \VV$ is the \emph{edge set}.  An edge from $i$ to $j$,
denoted by $(i,j)$, means that agent $j$ can send information to agent
$i$. For an edge $(i,j) \in\EE$, $i$ is called an \emph{in-neighbor}
of $j$, and $j$ is called an \emph{out-neighbor} of $i$.
A digraph $\GG'=(\VV,\EE')$ is a \emph{spanning subgraph} of a digraph
$\GG=(\VV,\EE)$ if $\EE'\subset\EE$.  A graph is \emph{undirected} if
$(i,j) \in \EE$ anytime $(j,i)\in\EE$.  Given digraphs
$\GG_i=(\VV,\EE_i)$, $i \in \until{m}$, defined on same node set, the
\emph{joint digraph} of these digraphs is the union $\cup_{i=1}^n
\GG_i=(\VV,\EE_1\cup\EE_2\cup \dots\cup\EE_m)$.  A \emph{directed
  path} is an ordered sequence of vertices such that any ordered pair
of vertices appearing consecutively is an edge of the digraph. A
\emph{directed tree} is an acyclic digraph with the following
property: there exists a node, called the root, such that any other
node of the digraph can be reached by one and only one directed path
starting at the root. A \emph{directed spanning tree} of a digraph is
a spanning subgraph that is a directed tree.  A digraph is called
\emph{strongly connected} if for every pair of vertices there is a
directed path between them.

A \emph{weighted digraph} is a triplet $\GG = (\VV
,\EE,\vect{\mathsf{A}})$, where $(\VV ,\EE )$ is a digraph and
$\vect{\mathsf{A}}\in\real^{N\times N}$ is a weighted
\emph{adjacency} matrix with the property that $a_{ij} >0$ if $(i, j)
\in\EE$ and $a_{ij} = 0$, otherwise. We use
$\Gamma(\vect{\mathsf{A}})$ to denote a digraph induced by a given
adjacency matrix $\vect{\mathsf{A}}$. A weighted digraph is
\emph{undirected} if $a_{ij} = a_{ji}$ for all $i,j\in\VV$.  The
\emph{weighted out-degree} and \emph{weighted in-degree} of a node
$i$, are respectively, $\text{d}^{\text{in}}(i ) =\sum^N_{j =1}a_{ji}$
and $\text{d}^{\text{out}} (i) =\sum^N_{j =1} a_{ij}$. We let
$\text{d}_{\max}^{\text{out}} = \underset{i \in \until{N}}{\max}
\text{d}^{\text{out}} (i)$ 
denote the
maximum  
weighted out-degree.
A digraph is
\emph{weight-balanced} if at each node $i\in\VV$, the weighted
out-degree and weighted in-degree coincide (although they might be
different across different nodes).  
The out-degree matrix $\vect{\mathsf{D}}^{\text{out}}$ is the diagonal matrix
with entries $\vect{\mathsf{D}}^{\text{out}}_{ii} = \text{d}^{\text{out}}(i)$,
for all $i\in\VV$. The \emph{(out-) Laplacian} matrix is $\vect{\mathsf{L}} =
\vect{\mathsf{D}}^{\text{out}} -\vect{\mathsf{A}}$. Note that
$\vect{\mathsf{L}}\vect{1}_N=0$. A weighted digraph $\GG$ is weight-balanced if
and only if $\vect{1}_N^T\vect{\mathsf{L}}=0$. Based on the structure of
$\vect{\mathsf{L}}$, at least one of the eigenvalues of $\vect{\mathsf{L}}$ is zero and
the rest of them have nonnegative real parts. 
We denote the eigenvalues of $\vect{\mathsf{L}}$ by $\lambda_i$, $i \in
\until{N}$, where $\lambda_1=0$ and
$\Re(\lambda_i)\leq\Re(\lambda_j)$, for $i<j$. For a strongly
connected digraph, zero is a simple eigenvalue of $\vect{\mathsf{L}}$. We
denote the eigenvalues of $\Sym{\vect{\mathsf{L}}}$ by $\hat{\lambda}_i$, $i
\in \until{N}$. For a strongly connected and weight-balanced digraph,
zero is a simple eigenvalue of $\Sym{\vect{\mathsf{L}}}$. For such a digraph,
we order the eigenvalues of $\Sym{\vect{\mathsf{L}}}$ as $\hat{\lambda}_1 =
0<\hat{\lambda}_2\leq\hat{\lambda}_3\leq\dots\leq\hat{\lambda}_N$.

\vspace{-6pt}
\subsection{Time-varying interactions via switched
  systems}\label{sec::SwtchNetModel}\vspace{-2pt}
  
Here, we introduce our model of networks with fixed number of agents
but time-varying interaction topologies. Let
$(\VV,\EE(t),\vect{\mathsf{A}}(t))$ be a time-varying digraph, where
the nonzero entries of the adjacency matrix are uniformly lower and
upper bounded, i.e., $a_{ij}(t)\in[\underline{a},\bar{a}]$, where
$0<\underline{a}\leq\bar{a}$, if $(j,i)\in\EE(t)$, and $a_{ij}=0$
otherwise. Our model of time-varying networks is then
$\GG(t)=\Gamma(\vect{\mathsf{A}}_{\sigma(t)})$, $t\geq 0$, with
$\sigma:[0,\infty)\to \mathcal{P} = \until{m}$ a piecewise constant
signal belonging to some switching set $\mathcal{S}$.  Here, $m$ can
be infinity. In our developments later, we provide precise
specifications for $\mathcal{S}$.  By piecewise constant, we mean a
signal that only has a finite number of discontinuities in any finite
time interval and that is constant between consecutive discontinuities
(no chattering). 
Without loss of generality, we assume that switching signals are
continuous from the right. The uniform stability of switched linear
systems with time-dependent switching signals (where uniformity refers
to the multiple solutions that can be obtained as the switching signal
ranges over a switching set) is characterized by the following result.
\begin{lem}[Asymptotic stability of switched linear systems implies
  exponential stability~\cite{JH:04}]\label{lem::swchHspa}
  For linear switched systems with trajectory-independent switching,
  uniform asymptotic stability is equivalent to exponential
  stability.
\end{lem}

We end this section by introducing the following notations. Given a
time-varying digraph, we denote by ${\cup}_{t_1}^{t_2}\GG(t)$ the
joint digraph in the time interval $[t_1,t_2)$ where
$t_1<t_2<+\infty$. We say a time-varying graph $\GG(t)$ is
\emph{jointly strongly connected} over the time-interval $[t_1,t_2)$
if ${\cup}_{t_1}^{t_2}\GG(t)$ is strongly connected. The time instants
at which the switching signal $\sigma$ is discontinuous are called
\emph{switching times} and are denoted by $t_0,t_1,t_2,\cdots$, where
$t_0=0$. We use $\vect{\mathsf{L}}_{\sigma}$ to represent the out-Laplacian of
the digraph $\Gamma(\vect{\mathsf{A}}_{\sigma})$.

\vspace{-6pt}
\section{Problem statement}\label{sec::ProbDef}\vspace{-2pt}

We consider a network of $N$ agents with single-integrator dynamics
given by
\begin{equation}\label{eq::AgentSingInt}
  \dot{x}^i=c^i,~~~i \in \until{N},
\end{equation}
where $x^i\in\real$ is the \emph{agreement state} and $c^i\in\real$ is
the \emph{driving command} of agent $i$. The network interaction
topology is modeled by a weighted digraph $\GG$.  Agent $i \in
\until{N}$ has access to a time-varying input signal
$u^i:[0,\infty)\to\real$.  The problem we are interested in solving is
the following.

\begin{prob}[Dynamic average consensus]\label{prob::main}
  Let $\GG$ be strongly connected and weight-balanced.  Design a
  distributed algorithm such that each agent's state $x^i(t)$
  asymptotically tracks the average $\avrg{u^j(t)}$ of the
  inputs.\boxend
\end{prob}

This problem finds numerous applications in networks of multiple
agents that have access to partial and evolving information, and aim
to combine it in a dynamic fashion.  Examples are numerous and include
data fusion, spatial estimation, and localization and mapping, to name
a few.  The algorithm design amounts to specifying a suitable driving
command $c^i$ for each agent $i \in \until{N}$. By distributed, we
mean that agent $i$ only interacts with its out-neighbors.  In
addition, we also consider variations of the problem above that are
intended to satisfy some practical issues that arise in using the
consensus algorithm in applications where the agent state corresponds
to a physical quantity such as position or velocity in motion
coordination of autonomous mobile agents.  In such applications, a
genuine concern is whether the command $c^i$ dictated by the consensus
algorithm can be implemented given the physical limitation of the
actuation systems. This motivates us to formulate the following
variation of Problem~\ref{prob::main}.

\begin{prob}[Dynamic average consensus with controllable rate of
  convergence]\label{prob::rate}
  Solve Problem~\ref{prob::main} such that each agent converges at its
  own desired rate of convergence.\boxend
\end{prob}

By giving a freedom to choose their desired rate of convergence, we
allow agents with limited control authority to opt for a slow rate of
convergence. We can also use the control over the individual rate of
convergence of agents in scheduling different time of arrivals. This
can benefit applications such as payload delivery or aerial
surveillance.  Although reducing the rate of convergence helps with
cases where control authority is limited, there is no guarantee that
control bounds, if present, would be satisfied.  This motivates us to
formulate the next problem.

\begin{prob}[Dynamic average consensus with limited control
  authority]\label{prob::sat}
  Solve Problem~\ref{prob::main} under bounded driving commands, i.e.,
  $\dot{x}^i=\sat{\bar{c}^i}{c^i}$ for all $i\in\until{N}$.\boxend
\end{prob}

Finally, we consider the problem of dynamic average consensus with
privacy preservation in the presence of adversaries.  Our motivation
to study such properties stems from the fact that privacy guarantees
on a distributed algorithm facilitate the agent participation in the
completion of cooperative tasks.  In an average consensus problem, the
privacy concern of agents can be local (e.g., some or all of the
agents do not want to reveal their local inputs to the outside world)
or global (e.g., all agents do not want to reveal their agreement
value to agents outside network).  We consider adversaries inside or
outside the network that do not interfere with the algorithm
implementation but seek to steal information about the inputs,
agreement value, or the agreement state trajectories of the individual
agents. The information these adversaries can access includes the time
history of intra network communication messages, partial or full
knowledge about the communication topology, and the algorithm design
parameters, and/or its initial conditions.

\begin{prob}[Dynamic average consensus with privacy
  preservation]\label{prob::priv}
  Solve Problems~\ref{prob::main}-\ref{prob::sat} such that the
  following privacy requirements are satisfied
  \begin{itemize}
  \item[(a)] the local inputs of the agents should not be revealed or
    be reconstructible by any adversary;
  \item[(b)] the agreement value should not be revealed to or be
    reconstructible by external adversaries;
  \item[(c)] the agreement state should not be revealed to or be
    reconstructible by any adversary. \boxend
   \end{itemize}
  \end{prob}

  For vector-valued inputs, one can apply the solution of
  Problems~\ref{prob::main}-\ref{prob::priv} in each dimension.

\vspace{-6pt}
\section{Dynamic average
  consensus}\label{sec::DyConsensus}\vspace{-2pt}

In this section, we introduce a distributed dynamic average consensus
algorithm that solves Problem~\ref{prob::main} with a steady-state
error for arbitrary time-varying input signals. We show that the
size of this error can be controlled using a design parameter and
that, for special classes of inputs, the steady-state error is zero. We also analyze the asymptotic correctness of the algorithm under time-varying interaction topologies and characterize the requirements on the stepsize for~discrete-time~implementations. 

\vspace{-6pt}
\subsection{Fixed interaction topology}\vspace{-2pt}

Here, we assume that the interaction topology of the network is fixed.
We propose the following distributed algorithm as our solution for
Problem~\ref{prob::main}
\begin{subequations}\label{eq::DC1}
  \begin{align}
    &\dot{x}^i = \dot{u}^i-\alpha
    (x^i-u^i)-\beta\SUM{j=1}{N}\vect{\mathsf{L}}_{ij}x^j-v^i,\label{eq::DC1-a}
    \\
    &\dot{v}^i=\alpha\beta\SUM{j=1}{N}\vect{\mathsf{L}}_{ij}x^j,\label{eq::DC1-b}
  \end{align}
\end{subequations}
where for $i \in \until{N}$, $x^i,v^i\in\real$ are
variables associated with agent~$i$. Also, $\vect{\mathsf{L}}$ is the
Laplacian of the digraph $\GG$ modeling the interaction topology.
This algorithm uses the last two terms of~\eqref{eq::DC1-a} as a
proportional integral feedback to impose agreement among neighboring
agents while these agents, because of the first two terms
of~\eqref{eq::DC1-a}, are moving towards their respective input
signal. Under suitable conditions on the communication topology,
explained below, this scheme results in each agents eventually
following the average of all the inputs across the network. The
constants $\alpha$, $\beta \in \real$ are design parameters that can
be used to tune the algorithm performance.  In the following, we study
the convergence and stability properties by using the equivalent
compact form below
\begin{subequations}\label{eq::DC1_shf}
  \begin{align}
    &\dvect{y} =
    -\alpha\vect{y}-\beta\vect{\mathsf{L}}\vect{y}-\vect{w},\label{eq::DC1_shf-a}
    \\
    &\dvect{w} =
    \alpha\beta\vect{\mathsf{L}}\vect{y}-\vect{\Pi}_N(\ddvect{u}+\alpha\dvect{u}).
    \label{eq::DC1_shf-b}
  \end{align}
\end{subequations}
where 
\begin{subequations}\label{eq::xvToyw}
  \begin{align}
    y^i &=x^i-\avrg{u^j}, \quad i \in \until{N},\label{eq::xTOy}\\
    \vect{w}&=\vect{v}-\bar{\vect{v}}, \quad
    \bar{\vect{v}}=\vect{\Pi}_N(\dvect{u}+\alpha\vect{u}).\label{eq::wTOv}
  \end{align}
\end{subequations}
Recall from Section~\ref{sec::ProbDef} that $x^i$ is the
agreement state of agent $i$. Thus, with the change of
variables~\eqref{eq::xTOy} we are transferring the desired equilibrium
of the system, in agreement state, to zero. We start our study by analyzing the stability and convergence
properties of the zero-system of~\eqref{eq::DC1_shf}, i.e.,
\begin{align}\label{eq::UF1C}
  \begin{bmatrix}
    \dvect{y}
    \\
    \dvect{w}
  \end{bmatrix} 
  = \vect{A}
  \begin{bmatrix}
    \vect{y}
    \\
    \vect{w}
  \end{bmatrix}, \quad \text{where } 
  \vect{A}=\begin{bmatrix}-\alpha\vect{I}_N-\beta\vect{\mathsf{L}}&-\vect{I}_N\\
    \alpha\beta\vect{\mathsf{L}}&\vect{0}
  \end{bmatrix}.
\end{align}
In the following, we show that the dynamical system~\eqref{eq::UF1C},
over a strongly connected and weight-balanced digraph, is stable and
convergent.

\begin{lem}[Asymptotic convergence
  of~\eqref{eq::UF1C}]\label{lem::UF1}
  Let $\GG$ be strongly connected and weight-balanced. For any $\alpha$,
  $\beta>0$, the trajectory of~\eqref{eq::UF1C} over~$\GG$ starting
  from any initial condition $\vect{y}(0), \vect{w}(0) \in\real^N$
  satisfies, 
  \begin{equation}\label{eq::UF1_equilibPoint}
    y^i(t)\to -\frac{\alpha^{-1}}{N}\SUM{j=1}{N}w^j(0),  \quad w^i(t)\to
    \avrg{w^j(0)}, \quad \text{as } t\to\infty, \quad \forall i \in \until{N},
  \end{equation}
  exponentially fast  with a rate of convergence upper
  bounded by $\min \{\alpha,\beta \Re(\lambda_2) \}$.
\end{lem}
\begin{pro}
  Consider the following change of variables where $\vect{r}=\frac{1}{\sqrt{N}}\vect{1}_N$ and $\vect{R}$ is such
  that $\vect{r}^\top\vect{R}=0$ and
  $\vect{R}^\top\vect{R}=\vect{I}_{N-1}$,
  \begin{equation}\label{eq::ChVarSim}
    \begin{bmatrix}
      \vect{p}
      \\
      \vect{q}
    \end{bmatrix}
    =
    \vect{T}_1 \vect{T}_2
    \begin{bmatrix}
      \vect{y}
      \\
      \vect{w}
    \end{bmatrix}, \quad \vect{T}_1=
    \begin{bmatrix}
      \vect{I}_N & \vect{0}
      \\
      \alpha \vect{I}_N & \vect{I}_N
    \end{bmatrix}, \quad
    \vect{T}_2
    =
    \begin{bmatrix}
      \vect{T}_3^\top&\vect{0}
      \\
      \vect{0}&\vect{T}_3^\top
    \end{bmatrix}, \quad \vect{T}_3 =
    \begin{bmatrix}
      \vect{r} & \vect{R}
    \end{bmatrix}.
  \end{equation}
   We partition the new
  variables as $\vect{p}=(p_1, \vect{p}_{2:N})$ and $\vect{q}=(q_1,
  \vect{q}_{2:N})$, where $p_1,q_1\in\real$ and
  $\vect{p}_{2:N},\vect{q}_{2:N}\in\real^{N-1}$.  Using~\eqref{eq::ChVarSim} the dynamics~\eqref{eq::UF1C} can be stated in the
  following equivalent form
  \begin{subequations}\label{eq::UF1C_Smlr}
    \begin{alignat}{2}
      \begin{bmatrix}
        \dot{p}_1
        \\
        \dot{q}_1
      \end{bmatrix}
      & = \Tvect{A}
      \begin{bmatrix}
        p_1
        \\
        q_1
      \end{bmatrix}, & \quad \Tvect{A}
        &=
      \begin{bmatrix}
        0&-1
        \\
        0&
        -\alpha
      \end{bmatrix}, \label{eq::UF1C_Smlr-a}
      \\
      \begin{bmatrix}
        \dvect{p}_{2:N}
        \\
        \dvect{q}_{2:N}
      \end{bmatrix}
      &= \Bvect{A}
      \begin{bmatrix}
        \vect{p}_{2:N}
        \\
        \vect{q}_{2:N}
      \end{bmatrix}, & \quad \Bvect{A} &=
      \begin{bmatrix}
        -\beta\vect{R}^\top\vect{\mathsf{L}}\vect{R}&-\vect{I}_{N-1}
        \\
        \vect{0}&-\alpha\vect{I}_{N-1}
      \end{bmatrix}.
      \label{eq::UF1C_Smlr-b}
    \end{alignat}
  \end{subequations}
  The eigenvalues of $\Tvect{A}$ are $0$ and $-\alpha$. The
  eigenvalues of the matrix $\Bvect{A}$ are $-\alpha$, with
  multiplicity $N-1$, and $-\beta\lambda_i$, with $i \in
  \{2,\dots,N\}$.  Recall that $\lambda_i$'s are eigenvalues of
  $\vect{\mathsf{L}}$. For a strongly connected digraph, $\lambda_1=0$ and the
  rest of the eigenvalues have positive real parts.  Therefore, for
  $\alpha, \beta>0$, the dynamical system~\eqref{eq::UF1C_Smlr}, and
  equivalently~\eqref{eq::UF1C}, is a stable linear system.

  The null-space of the system matrix $\vect{A}$ is spanned by
  $({\vect{1}_N},-\alpha{\vect{1}_N})$, the eigenvector associated
  with zero eigenvalue. Therefore,~\eqref{eq::UF1C} converges
  exponentially fast to the set
  \begin{equation}\label{eq::UF1_equilib}
    \{(\vect{y},\vect{w})\,|\,\vect{y}=\mu\vect{1}_N,~\vect{w}
    = -\mu\alpha\vect{1}_N,\quad\mu\in\real\}. 
  \end{equation}
  Left multiplying both sides of~\eqref{eq::UF1C} by
  $\Diag{{\vect{0}_N}^\top,{\vect{1}_N}^\top}$ and invoking the
  weight-balanced property of the digraph, we obtain
  $ \SUM{i=1}{N}\dot{w}^i =0$, and therefore,
  \begin{align}\label{eq::UF1_equilibSum} 
    &\SUM{i=1}{N}w^i(t)=\SUM{i=1}{N}w^i(0),~~\forall~t\geq 0.
  \end{align}
  The combination of~\eqref{eq::UF1_equilib}
  and~\eqref{eq::UF1_equilibSum} yields that, from any initial condition
  $\vect{y}(0), \vect{w}(0) \in\real^N$, the trajectory of the
  dynamical system~\eqref{eq::UF1C}
  satisfies~\eqref{eq::UF1_equilibPoint}, exponentially fast.  Based
  on~\eqref{eq::rate_of_conv}, the rate of convergence is $\min \{
  \alpha,\beta \Re(\lambda_2) \}$.
\end{pro}

The next result further probes into the properties of the dynamical
system~\eqref{eq::UF1C} by upper bounding the difference between the
state $y^i$ of agent $i$ at any time $t$ and the equilibrium
value. This bound is instrumental later in the characterization of the
steady-state error of~\eqref{eq::DC1}.

\begin{lem}[Upper bound on trajectories
  of~\eqref{eq::UF1C}]\label{eq::UF1_wyBound}
  Under the assumptions of Lemma~\ref{lem::UF1}, the following bound
  holds for each $i \in \until{N}$,
  \begin{equation*}
    \left| y^i(t)+\frac{\alpha^{-1}}{N}\sum_{j=1}^Nw^j(0)\right|\leq
    \Lnorm\vect{y}(t) + \alpha^{-1}\vect{r}\vect{r}^\top
    \vect{w}(0)\Rnorm\leq s(t), 
  \end{equation*}
  where
  \begin{align}\label{eq::y_upper_bound}
    s(t) &= (\e{-\alpha t} +\e{-\beta\Hlambda_2
      t})\Lnorm\vect{y}(0)\Rnorm+\alpha^{-1}\e{-\alpha
      t}\Lnorm\vect{w}(0)\Rnorm\nonumber
    \\
    & \quad +
    \begin{cases}
      (\beta\Hlambda_2-\alpha)^{-1}(\e{-\alpha
        t}-\e{-\beta\Hlambda_2t})\left(\alpha\Lnorm \vect{y}(0)\Rnorm
        + \Lnorm\vect{w}(0)\Rnorm\right) ,& \text{if }
      \alpha\neq\beta\Hlambda_2,
      \\
      t\e{-\beta\Hlambda_2t}\left(\alpha\Lnorm\vect{y}(0)\Rnorm +
        \Lnorm\vect{w}(0)\Rnorm\right), & \text{if } \alpha =
      \beta\Hlambda_2.
    \end{cases}
  \end{align}
 \end{lem}
\begin{pro}
  The solution of the state equation~\eqref{eq::UF1C_Smlr} from any
  initial condition $\vect{y}(0), \vect{w}(0) \in\real^N$ is $(
  p_1(t), q_1(t), \vect{p}_{2:N}(t), \vect{q}_{2:N}(t)) =
  \vect{\Omega}(t) ( p_1(0), q_1(0), \vect{p}_{2:N}(0),
  \vect{q}_{2:N}(0))$, where
  \begin{align}\label{eq::sol_UF1C_Smlr}
    \vect{\Omega}(t) & =
    \begin{bmatrix}
      1&\alpha^{-1} (\e{-\alpha t}-1 )&0&0
      \\
      0&\e{-\alpha t}&0&0
      \\
      \vect{0}&\vect{0}&\vect{\Phi}(t,0)&-\int_{0}^t
      \vect{\Phi}(t,\tau)\e{-\alpha\tau}d\tau
      \\
      \vect{0}&\vect{0}&\vect{0}&\e{-\alpha t}\vect{I}_{N-1}
    \end{bmatrix},
  \end{align}
  and $\vect{\Phi}(t,\tau) =
  \e{-\beta\vect{R}^\top\vect{\mathsf{L}}\vect{R}(t-\tau)}$. Now,
  from~\cite[Fact 11.15.7, item xvii]{DSB:09}, we deduce
  \begin{equation}\label{eq::R'LR_bund}
    \Lnorm\vect{\Phi}(t,\tau)\Rnorm=\Lnorm\e{-\beta\vect{R}^\top\vect{\mathsf{L}}
      \vect{R}(t-\tau)}\Rnorm\leq 
    \e{-\beta\Hlambda_2 (t-\tau)},
  \end{equation}
  and hence
  \begin{equation}\label{eq::EXP12}
    \Lnorm\int_{0}^t\vect{\Phi}(t,\tau)\e{-\alpha\tau}d\tau\Rnorm\leq
    \int_{0}^t\e{-\beta\Hlambda_2(t-\tau)}\e{-\alpha\tau}d\tau.
  \end{equation}
  Now, using the change of variables~\eqref{eq::ChVarSim}, one has
  \begin{equation}\label{eq::main_sol}
    \vect{y}(t)=\vect{S}_{11}\vect{y}(0)+\vect{S}_{12}\vect{w}(0),
  \end{equation}
  where
  \begin{subequations}\label{eq::S}
    \begin{align}
      &\vect{S}_{11}=\e{-\alpha
        t}\vect{r}\vect{r}^\top+\vect{R}\vect{\Phi}(t,0)\vect{R}^\top-\alpha
      \vect{R}(\int_{0}^t\vect{\Phi}(t,\tau)\e{-\alpha\tau}d\tau)
      \vect{R}^\top,
      \\
      &\vect{S}_{12}=(-\alpha^{-1}+\alpha^{-1}\e{-\alpha
        t})\vect{r}\vect{r}^\top-\vect{R}(\int_{0}^t\vect{\Phi}(t,\tau)
      \e{-\alpha\tau}d\tau)\vect{R}^\top.
    \end{align}    
  \end{subequations}
  The result now follows from using~\eqref{eq::R'LR_bund}
  and~\eqref{eq::EXP12} to bound the expression~\eqref{eq::main_sol}.
\end{pro}

Next, using the results guaranteed by Lemma~\ref{eq::UF1_wyBound}
we study the convergence and stability properties of our proposed
dynamic average consensus algorithm~\eqref{eq::DC1}.  We start by
establishing an upper bound on its tracking error for any given initial
condition.

\begin{thm}[Upper bound on the tracking error
  of~\eqref{eq::DC1}]\label{thm::DC1_bound}
  Let $\GG$ be strongly connected and weight-balanced. Each agent has
  a piecewise continuously differentiable input $u^i(t)$. For
  $\alpha$, $\beta>0$, the trajectory of the algorithm~\eqref{eq::DC1}
  over~$\GG$ starting from any initial condition $\vect{x}(0),
  \vect{v}(0) \in\real^N$ satisfies, for all $i \in \until{N}$,
  \begin{align}\label{eq::yDC_bound}
    \left| x^i(t)-\avrg{u^j(t)} +
      \frac{\alpha^{-1}}{N}\sum_{j=1}^Nv^j(0)\right|\leq&s(t) +
    \int_0^t
    \e{-\beta\Hlambda_2(t-\tau)}\Lnorm\vect{\Pi}_N\dvect{u}(\tau)\Rnorm
    d\tau+\\
    &
    \begin{cases}
      (\beta\Hlambda_2-\alpha)^{-1}(\e{-\alpha
        t}-\e{-\beta\Hlambda_2t}) \Lnorm\dvect{u}(0)\Rnorm , & \text{
        if } \alpha\neq\beta\Hlambda_2,
      \\
      t\e{-\beta\Hlambda_2t} \Lnorm\dvect{u}(0)\Rnorm , & \text{ if }
      \alpha=\beta\Hlambda_2,
    \end{cases}\nonumber
  \end{align}
  where $s(t)$ is defined in~\eqref{eq::y_upper_bound}, and $\vect{y}$
  and $\vect{w}$ are defined in~\eqref{eq::xvToyw}.
\end{thm}
\begin{pro}
  Using the change of
  the variables~\eqref{eq::ChVarSim} we can
  represent~\eqref{eq::DC1_shf}, an equivalent representation of~\eqref{eq::DC1}, in the following equivalent form where $\Tvect{A}$ and $\Bvect{A}$ are defined
  in~\eqref{eq::UF1C_Smlr},
  \begin{subequations}\label{eq::DC1_Smlr}
    \begin{align}
      &
      \begin{bmatrix}
        \dot{p}_1\\
        \dot{q}_1
      \end{bmatrix}=\Tvect{A}\begin{bmatrix}
        p_1\\
        q_1
      \end{bmatrix},
          \\
      &
      \begin{bmatrix}\dvect{p}_{2:N}\\
        \dvect{q}_{2:N}\end{bmatrix}
      =\Bvect{A}\begin{bmatrix}\vect{p}_{2:N}\\
        \vect{q}_{2:N}\end{bmatrix}-\begin{bmatrix}\vect{0}\\
        \vect{R}^\top\end{bmatrix}(\ddvect{u} +
      \alpha\dvect{u}),\label{eq::DC1_Smlr-b}
    \end{align}
  \end{subequations}
  For any given initial conditions, the
  solution of the state equation~\eqref{eq::DC1_Smlr} is
  \begin{equation*}
    \begin{bmatrix}
      p_1(t)\\
      q_1(t)\\
      \vect{p}_{2:N}(t)\\
      \vect{q}_{2:N}(t)
    \end{bmatrix}\!\!=\!\vect{\Omega}(t)\!\!\begin{bmatrix}
      p_1(0)\\
      q_1(0)\\
      \vect{p}_{2:N}(0)\\
      \vect{q}_{2:N}(0)
    \end{bmatrix}\!\!-\!\!\begin{bmatrix}0\\0\\
      \int_{0}^t\vect{\Phi}(t,\tau)\e{-\alpha\tau}d\tau\,(\vect{q}_{2:N}(0)
      \!\!+\!\!\vect{R}^\top\dvect{u}(0))\!-\!\int_{0}^t\vect{\Phi}(t,\tau)
      \vect{R}^\top\dvect{u}(\tau)d\tau\\
      -\vect{R}^\top\dvect{u}(0)+\vect{R}^\top\dvect{u}(t)
    \end{bmatrix}\!,
  \end{equation*}
  where $\vect{\Omega}(t)$ is defined
  in~\eqref{eq::sol_UF1C_Smlr}. Recalling the change of
  variables~\eqref{eq::ChVarSim}, we have
  \begin{equation}\label{eq::SolDC1_y}
    \vect{y}(t) =
    \vect{S}_{11}\vect{y}(0)+\vect{S}_{12}\vect{w}(0) -
    \vect{R}\int_{0}^t\vect{\Phi}(t,\tau) 
    \e{-\alpha\tau}d\tau\,\vect{R}^\top\dvect{u}(0) +
    \vect{R}\int_{0}^t\vect{\Phi}(t,\tau)  
    \vect{R}^\top\dvect{u}(\tau)d\tau, 
  \end{equation}
  where $\vect{S}_{11}$ and $\vect{S}_{12}$ are defined
  in~\eqref{eq::S}. Note that~\eqref{eq::wTOv} implies that
  $\SUM{i=1}{N}w^i(0)=\SUM{i=1}{N}v^i(0)$. Notice also that
  $\vect{R}^\top = \vect{R}^\top\vect{\Pi}_N$, and $\Lnorm
  \vect{R}\Rnorm = \Lnorm \vect{R}^\top\Rnorm =
  \sigma_{\max}(\vect{R})=1$. Then, by
  recalling~\eqref{eq::R'LR_bund}, it is straightforward to show
  that~\eqref{eq::yDC_bound} is satisfied.
\end{pro}

The next result shows that, for input signals whose orthogonal
projection into the agreement space are essentially bounded, the
algorithm~\eqref{eq::DC1} solves Problem~\ref{prob::main} with a
bounded steady-state error.

\begin{cor}[The algorithm~\eqref{eq::DC1} solves
  Problem~\ref{prob::main}]\label{cor::DC1Prob1Sol}
  Let $\GG$ be strongly connected and weight-balanced.  Assume that
  the derivatives of the inputs of the network satisfy
  $\|\vect{\Pi}_N\dvect{u}\|_{\text{ess}} = \gamma<\infty$. Then, for
  any $\alpha,\beta>0$ the algorithm~\eqref{eq::DC1}
  over~$\GG$ initialized at $x^i(0),v^i(0)\in \real$
  such that $\SUM{i=1}{N}v^i(0)=0$ solves Problem~\ref{prob::main}
  with an upper-bounded steady-state error. Specifically,
  \begin{equation}\label{eq::DC1_ultimate_bound} 
    \lim_
    {t\to\infty}\sup\left|
      x^i(t)-\avrg{u^j(t)}\right|\leq (\beta\Hlambda_2)^{-1}\gamma,
    \quad i \in \until{N}. 
  \end{equation}
\end{cor}
\begin{pro}
  In Theorem~\ref{thm::DC1_bound}, for a strongly connected and
  weight-balanced digraph, we showed that the trajectories of the
  algorithm~\eqref{eq::DC1}, for any $x^i(0),v^i(0)\in\real$, $i \in
  \until{N}$, satisfy the bound~\eqref{eq::yDC_bound}.   Then, we can easily deduce~\eqref{eq::DC1_ultimate_bound}
  from~\eqref{eq::yDC_bound} using
  \begin{equation*}
    \int_0^t
    \e{-\beta\Hlambda_2(t-\tau)}
    \Lnorm\vect{\Pi}_N \dvect{u}(\tau)\Rnorm d\tau\leq (\beta\Hlambda_2)^{-1}
    (1-\e{-\beta\Hlambda_2 
      t})\gamma.
  \end{equation*}

\end{pro}

\begin{rem}[Effect of faulty initial conditions]\label{rem::v(0)}
  The condition $\SUM{i=1}{N}v^i(0)=0$ of
  Corollary~\ref{cor::DC1Prob1Sol} can be easily satisfied if each
  agent starts at $v^i(0)=0$. This is a mild requirement because $v^i$
  is an internal state for agent $i$, and therefore it is not affected
  by imperfect communication errors. Additionally, for large networks,
  if we assume that the initialization error is zero-mean Gaussian
  noise, we can expect $\SUM{i=1}{N}v^i(0)=0$.
  \boxend
\end{rem}

\begin{rem}[Tuning the performance of~\eqref{eq::DC1} via design
  parameters]\label{re:performance-regulation}
  Corollary~\ref{cor::DC1Prob1Sol} shows that to reduce the nonzero
  steady-state error, one can either increase the graph connectivity
  (larger $\Hlambda_2$) or use a larger value of $\beta$.  The
  parameter $\alpha$ can also be exploited to regulate the algorithm
  performance.  According to the bound~\eqref{eq::yDC_bound} the rate
  of convergence of the transient behavior is governed by
  $\min\{\alpha,\beta\Hlambda_2\}$. If one is forced to use large
  $\beta\Hlambda_2$ to reduce the steady-state error, then $\alpha$
  can fulfill the role of regulating the rate of convergence of the
  algorithm.
  \boxend
\end{rem}

\begin{rem}[Comparison with input requirements of the solutions in the
  literature] 
  In order to guarantee bounded steady-state tracking error, the solution we offer for Problem~\ref{prob::main}
  through Corollary~\ref{cor::DC1Prob1Sol} only requires that the
  projection of the network's aggregated input derivative vector into
  the agreement space is bounded. This is more general than the
  requirements in the literature, which generally ask for bounded
  input and/or bounded derivatives
  (e.g.,~\cite{ROS-JSS:05,RAF-PY-KML:06,MZ-SM:08a}).  \boxend
\end{rem}

In the following, we identify conditions involving the inputs and their
derivatives under which the algorithm~\eqref{eq::DC1}
solves Problem~\ref{prob::main} with zero steady-state error.

\begin{lem}[Conditions on inputs for zero steady-state error
  of~\eqref{eq::DC1}]\label{eq::Cnvrg}
  Let $\GG$ be strongly connected and weight-balanced. Assume there
  exists $\alpha >0$ such that, for all $i \in \until{N}$, one of the
  following conditions are satisfied
  \begin{itemize}
  \item[(a)] $\dot{u}^i(t)+\alpha u^i(t)$ converges to a common
    function $l(t)$ as $t\to \infty$;
  \item[(b)] $\ddot{u}^i(t)+\alpha \dot{u}^i(t)$ converges to a common
    function $l(t)$ as $t\to \infty$.
  \end{itemize}
  Then, the algorithm~\eqref{eq::DC1} over~$\GG$ with the given
  $\alpha$, and $x^i(0),v^i(0)\in\real$ such that
  $\SUM{i=1}{N}v^i(0)=0$, for any $\beta>0$, makes $x^i(t)\to
  \avrg{u^j(t)}$ as $t\to\infty$, for all $i \in \until{N}$.
\end{lem}
 \begin{pro}
Using the change of variables~\eqref{eq::xTOy} we can
  represent~\eqref{eq::DC1} in the following equivalent compact form
  \begin{subequations}\label{eq::DC1_shfx}
    \begin{align}
      &\dvect{y} = -\alpha\vect{y}-\beta\vect{\mathsf{L}}\vect{y}-\vect{v}
      +\vect{\Pi}_N(\dvect{u}+\alpha\vect{u}), \label{eq::DC1_shfx-a}
      \\
      &\dvect{v}=\alpha\beta\vect{\mathsf{L}}\vect{y}. \label{eq::UD-CF1_shfx-b}
    \end{align}
  \end{subequations}
  When condition (a) holds we have
  $\vect{\Pi}_N(\dvect{u}+\alpha\vect{u})\to\vect{0}$, as
  $t\to\infty$.  Then,~\eqref{eq::DC1_shfx} is a linear system with a
  vanishing input $\vect{\Pi}_N(\dvect{u}+\alpha\vect{u})$. Therefore,
  it converges to the equilibrium of its zero-system. In light of
  Lemma~\ref{lem::UF1}, we conclude that $y^i(t)\to
  -\frac{\alpha^{-1}}{N}\SUM{j=1}{N}v^j(0)$ asymptotically for all $i
  \in \until{N}$. However, due to initialization requirement we have
  $\SUM{i=1}{N}v^i(0)=0$. As a result $x^i(t)\to \avrg{u^j(t)}$
  globally asymptotically for $i \in \until{N}$.
  
  When condition (b) holds we have
  $\vect{\Pi}_N(\ddvect{u}+\alpha\dvect{u})\to\vect{0}$, as
  $t\to\infty$.  Recall~\eqref{eq::DC1_shf} the equivalent
  representation of~\eqref{eq::DC1}.  It is a linear system with a
  vanishing input $\vect{\Pi}_N(\dvect{u}+\alpha\vect{u})$. Then,
  using a similar argument used for~\eqref{eq::DC1_shfx} above, we can
  show that in~\eqref{eq::DC1_shf} $y^i(t)\to
  -\frac{\alpha^{-1}}{N}\SUM{j=1}{N}w^j(0)$ asymptotically for all $i
  \in \until{N}$.  Using~\eqref{eq::wTOv}, we can show
  $\SUM{i=1}{N}w^i(0)=\SUM{i=1}{N}v^i(0)$. As a result $x^i(t)\to
  \avrg{u^j(t)}$ globally asymptotically for $i \in \until{N}$.
\end{pro}

\begin{rem}[Inputs that satisfy the conditions of
  Lemma~\ref{eq::Cnvrg}]
  The classes of inputs in Lemma~\ref{eq::Cnvrg} depend on the parameter
  $\alpha$ which must be known by each agent in order to obtain zero
  steady-state error. There are classes of inputs that satisfy the conditions
  regardless of the value of $\alpha$, such as static inputs and
  dynamic inputs which differ from one another by static values. For
  these classes of inputs,
  $\vect{\Pi}_N(\ddvect{u}+\alpha\dvect{u})=\vect{0}$, and the
  convergence is exponential with rate
  $\min\{\alpha,\beta\Re(\lambda_2)\}$. \boxend
\end{rem}

\vspace{-6pt}
\subsection{Time-varying interaction
  topologies}\vspace{-2pt}

In this section, we analyze the stability and convergence properties
of the dynamic average consensus algorithm~\eqref{eq::DC1} over
networks with changing interaction topology.  Changes can be due to
unreliable transmission, limited communication/sensing range, or
obstacles.  Let $(\VV,\EE(t),\vect{\mathsf{A}}(t))$ be a time-varying
digraph, where the nonzero entries of the adjacency matrix are
uniformly lower and upper bounded (i.e., $a_{ij}(t)\in[\underline{a},
\bar{a}]$, where $0<\underline{a}\leq\bar{a}$, if $(j,i)\in\EE(t)$,
and $a_{ij}=0$ otherwise).  Intuitively one can expect that consensus
in switching networks will occur if there is occasional enough flow of
information from every node in the network to every other node. Then,
according to Section~\ref{sec::SwtchNetModel}, in order to describe
our switching network model, we start by specifying the set of
admissible switching signals.

\begin{defn}[\emph{Admissible} switching set
  $\mathcal{S}_{\mathrm{admis}}$]
  An admissible switching set $\mathcal{S}_{\mathrm{admis}}$ is a set
  of piecewise constant switching signals $\sigma: [0,\infty)
  \rightarrow \mathcal{P}$ with some dwell time $t_L$ (i.e.,
  $t_{k+1}-t_{k}>t_L>0$, for all $k=0,1,\dots$) such that
  \begin{itemize}
  \item the induced digraph $\Gamma(\vect{\mathsf{A}}_{\sigma(t)})$
    is weight-balanced for $t\geq t_0$;
  \item the number of contiguous, nonempty, uniformly bounded
    time-intervals $[t_{i_j},t_{i_{j+1}})$, $j=1,2,\dots$, starting
    at $t_{i_1}=t_0$, with the property that
    ${\cup}_{t_{i_j}}^{t_{i_{j+1}}}\Gamma(\vect{\mathsf{A}}_{\sigma(t)})$
    is a jointly strongly connected digraph goes to infinity as
    $t\to\infty$. \boxend\end{itemize}
\end{defn}

Our model of network with switching topology is then
$\Gamma(\vect{\mathsf{A}}_\sigma)$, with $\sigma \in
\mathcal{S}_{\mathrm{admis}}$.
The algorithm~\eqref{eq::DC1}, after applying the change of
variables~\eqref{eq::xvToyw}, is represented in compact form as
follows
\begin{align}\label{eq::DC1Switch}
  \begin{bmatrix}
    \dvect{y}\\
    \dvect{w}
  \end{bmatrix}=\vect{A}_{\sigma(t)}\begin{bmatrix}
    \vect{y}\\
    \vect{w}
  \end{bmatrix}-\begin{bmatrix}\vect{0}
    \\
    \vect{\Pi}_N(\ddvect{u}+\alpha\dvect{u})
  \end{bmatrix}, \quad \vect{A}_{\sigma(t)}
  = \begin{bmatrix}-\alpha\vect{I}_N-\beta
    \vect{\mathsf{L}}_{\sigma(t)}&-\vect{I}_N
    \\
    \alpha\beta \vect{\mathsf{L}}_{\sigma(t)}&\vect{0}.
  \end{bmatrix}.
\end{align}
Similarly to our analysis of the algorithm over fixed interaction
topologies, we start by examining the zero-system
of~\eqref{eq::DC1Switch}, i.e.,
\begin{equation}\label{eq::UF1Switch}
 \begin{bmatrix}
   \dvect{y}\\
   \dvect{w}
  \end{bmatrix}=\vect{A}_{\sigma(t)}\begin{bmatrix}
    \vect{y}\\
    \vect{w}
  \end{bmatrix}.
\end{equation}
The following result analyzes the convergence and stability properties
of the switched dynamical system~\eqref{eq::UF1Switch} when the
switching signal $\sigma\in\mathcal{S}_{\mathrm{admis}}$.

\begin{lem}[Asymptotic convergence
  of~\eqref{eq::UF1Switch}]\label{lem::swtch_stable_UF}
  Let $\sigma\in \mathcal{S}_{\mathrm{admis}}$ and consider $\GG(t)=
  \Gamma(\vect{\mathsf{A}}_{\sigma(t)})$ for $t\geq0$.
  Then, for any $\alpha,\beta>0$, the trajectory of the
  algorithm~\eqref{eq::UF1Switch} starting from any initial condition
  $\vect{y}(0), \vect{w}(0) \in\real^N$
  satisfies~\eqref{eq::UF1_equilibPoint}, exponentially fast.
\end{lem}
\begin{pro}
  Using the change of the variables~\eqref{eq::ChVarSim}, we can
  represent~\eqref{eq::UF1Switch} in the equivalent form~\eqref{eq::UF1C_Smlr} in which $\Bvect{A}$ and $\vect{\mathsf{L}}$ are replaced by $\Bvect{A}_{\sigma(t)}$ and $\vect{\mathsf{L}}_{\sigma(t)}$, respectively.  
  We can write $\dvect{p}$ as follows
  \begin{equation}\label{eq::p_dot_swch}
    \dvect{p}=-\vect{T}_3^\top\vect{\mathsf{L}}_\sigma\vect{T}_3\vect{p}-\vect{q}.
  \end{equation}
  We can look at this dynamical equation as a linear system with input
  $\vect{q}$ which vanishes exponentially fast (notice that
  $\dvect{q}=-\alpha\vect{q}$). Next, we examine the stability of
  zero-system of~\eqref{eq::p_dot_swch}. Under the state
  transformation $\vect{\eta}=\vect{T}_3\vect{p}$, this zero-system
  can be represented in the following equivalent form
  \begin{equation}\label{eq::static_cons}
    \dvect{\eta}=-\vect{\mathsf{L}}_\sigma\vect{\eta}.
  \end{equation}
  According to \cite[Theorem 2.33]{WR-RWB:08}, when the switching
  signal $\sigma$ is such that the number of contiguous, nonempty,
  uniformly bounded time-intervals $[t_{i_j},t_{i_{j+1}})$,
  $j=1,2,\dots$, starting at $t_{i_1}=t_0$, with the property that
  ${\cup}_{t_{i_j}}^{t_{i_{j+1}}}\Gamma(\vect{\mathsf{A}}_{\sigma(t)})$
  has a spanning tree, then~\eqref{eq::static_cons} asymptotically
  achieves consensus. Invoking this result, we can conclude that for
  $\sigma \in \mathcal{S}_{\mathrm{admis}}$, the trajectories
  of~\eqref{eq::static_cons} converge asymptotically to
  $\avrg{\eta_j(0)}$ where $\eta_i(0)$ is the $i$th element of
  $\vect{\eta}(0)$. For zero-system of~\eqref{eq::p_dot_swch}, this is
  equivalent to $p_1(t)\to p_1(0)$ and $\vect{p}_{2:N}(t)\to \vect{0}$
  uniformly asymptotically for all
  $\sigma\in\mathcal{S}_{\mathrm{admis}}$. The switching signal
  $\sigma\in\mathcal{S}_{\mathrm{admis}}$ is a trajectory-independent
  (it is time-dependent) switching signal. Then,
  Lemma~\ref{lem::swchHspa} implies that the convergence of the zero
  system of~\eqref{eq::p_dot_swch} is indeed globally uniformly
  exponentially fast.  Using input-to-state stability results (see
  \cite{EDS-YW:95,LV-DC-DL:05}), then we can conclude that
  in~\eqref{eq::p_dot_swch}, $p_1(t)\to p_1(0)$ and
  $\vect{p}_{2:N}(t)\to \vect{0}$ as $t\to\infty$ uniformly globally
  exponentially.  Recall the change of variable~\eqref{eq::ChVarSim},
  then it is easy to show that for~\eqref{eq::UF1Switch} we also
  have~\eqref{eq::UF1_equilibPoint}.
\end{pro}

Obtaining an explicit value for the rate of convergence
of~\eqref{eq::UF1Switch} for all possible $\sigma \in
\mathcal{S}_{\mathrm{admis}}$ is not straightforward.  However, we can
show that the rate of convergence is upper bounded by
$\underset{p\in\PP}{\max}\left(\Re(\lambda_{p2})\right)$, where
$\lambda_{p2}$ is the eigenvalue of $\vect{\mathsf{L}}_p$ with smallest nonzero
real part.  The following result relates the upper bound on the
difference between the state $y^i(t)$ of agent $i$ at any time $t$ and
the final agreement value to the rate of convergence
of~\eqref{eq::UF1_equilibPoint}.

\begin{lem}[Upper bound on trajectories
  of~\eqref{eq::UF1Switch}]\label{lem::swtch_bound_UF}
  Under the assumptions of Lemma~\ref{lem::swtch_stable_UF}, the
  following bound holds for each $i\in\{1,\cdots,N\}$,
  \begin{align}\label{eq::ySwtch_bound}
    \left| y^i(t)+\frac{\alpha^{-1}}{N}\sum_{i=1}^Nw^i(0)
    \right|\leq\Lnorm\vectT{y}(t)+\alpha^{-1}\vect{r}
    \vect{r}^\top\vectT{w}(0)\Rnorm \leq \hat{s}(t),
  \end{align}
  where $\hat{s}(t)$ is the same as $s(t)$
  in~\eqref{eq::y_upper_bound} only $\Hlambda_2$ is replaced by
  $\Hlambda_{\sigma}>0$ where $\Hlambda_{\sigma}$ satisfies
  \begin{equation}\label{eq::Hlambda_sigma_Def}
    \Lnorm\e{-\beta
      \vect{R}^\top\vect{\mathsf{L}}_{\sigma(t)}\vect{R}(t-t_0)}\Rnorm\leq
    \kappa\e{-\beta\Hlambda_{\sigma}(t-t_0)},~~~\forall t\geq t_0\geq
    0,\end{equation}
  for some finite $0<\kappa$.
\end{lem}
\begin{pro}
  We follow the same steps of the proof of
  Lemma~\ref{eq::UF1_wyBound}. The only difference is that the norm
  bound~\eqref{eq::R'LR_bund} of the transition matrix of
  $\dvect{p}_{2:N}$ state equation has to be modified, as explained
  below.  We showed in the proof of Lemma~\ref{lem::swtch_bound_UF}
  that when $\sigma\in \mathcal{S}_{\mathrm{admis}}$ for all $t\geq
  t_0$, the zero-system of~\eqref{eq::p_dot_swch} is exponentially
  stable. Therefore, there exist positive $\Hlambda_{\sigma}$ and
  $\kappa$ such that
  \begin{align*}
    \Lnorm\Phi(t,t_0) = \e{-\beta
      \vect{R}^\top\vect{\mathsf{L}}_{\sigma(t)}\vect{R}(t-t_0)}\Rnorm\leq
    \kappa\e{-\beta\Hlambda_{\sigma}(t-t_0)},~~~\forall t\geq t_0\geq
    0.
  \end{align*}
  As a result, in the case of switched dynamical systems,
  in~\eqref{eq::EXP12} $\Hlambda_2$ is replaced by
  $\Hlambda_{\sigma}$. Then, from~\eqref{eq::main_sol} we can deduce
  the bound~\eqref{eq::ySwtch_bound}.
\end{pro}

In light of Lemma~\ref{lem::swtch_bound_UF}, the extension of the
results on the stability analysis and ultimate convergence error bound
of the algorithm~\eqref{eq::DC1} over fixed interaction topologies to
switching networks whose switching signal
$\sigma\in\mathcal{S}_{\text{admis}}$ is straightforward. For such
switching networks, Theorem~\ref{thm::DC1_bound} and
Corollary~\ref{cor::DC1Prob1Sol} are valid, with the only change of
replacing $\beta\Hlambda_2$ by $\beta\Hlambda_{\sigma}$,
cf.~\eqref{eq::Hlambda_sigma_Def}, in the statement. Because of
Lemma~\ref{lem::swtch_stable_UF}, the proof that Lemma~\ref{eq::Cnvrg}
applies to switched networks with
$\sigma\in\mathcal{S}_{\text{admis}}$ is straightforward. For the sake
of brevity the detailed statements and proofs are omitted.

\vspace{-6pt}
\subsection{Discrete-time implementation over fixed interaction
  topologies}\label{sec::Discrete}\vspace{-2pt}

Here, we study a discrete-time algorithm that solves
Problem~\ref{prob::main} with non-zero steady-state error. In doing
so, we are motivated by the aim of understanding the differences and
connections between continuous- and discrete-time systems for
multi-agent systems and by practical considerations regarding
algorithm implementability.  Given a stepsize $\delta>0$, for $i \in \until{N}$, consider
\begin{subequations}\label{eq::DCDisc}
  \begin{align}
    z^i(k+1) &= z^i(k)-\delta\alpha
    z^i(k)-\delta\beta\SUM{j=1}{N}\vect{\mathsf{L}}_{ij}(z^j(k)+u^j(k))-\delta
    v^i(k),\label{eq::DCDisc-a}
    \\
    v^i(k+1) &= v^i(k)+\delta\alpha\beta
    \SUM{j=1}{N}\vect{\mathsf{L}}_{ij}(z^j(k)+u^j(k)), \label{eq::DCDisc-b}
    \\
    x^i(k)&= z^i(k)+u^i(k). \label{eq::DCDisc-c}
  \end{align}
\end{subequations}
Using~\eqref{eq::DCDisc-c} to obtain
$z^i(k)=x^i(k)-u^i(k)$, and substituting this in~\eqref{eq::DCDisc-a}
and~\eqref{eq::DCDisc-b}, we obtain
\begin{subequations}\label{eq::DCDisc_smlr}
  \begin{align}
    &x^i(k+1) = x^i(k)-\delta\alpha
    (x^i(k)-u^i(k))-\delta\beta\SUM{j=1}{N}\vect{\mathsf{L}}_{ij}x^j(k)-\delta
    v^i(k)+\Delta u^i(k),\label{eq::DCDisc_smlr-a}
    \\
    &v^i(k+1) =
    v^i(k)+\delta\alpha\beta\SUM{j=1}{N}\vect{\mathsf{L}}_{ij}x^j(k), 
    \label{eq::DCDisc_smlr-b}
  \end{align}
\end{subequations}
where $\Delta u^i(k) =u^i(k+1)-u^i(k)$. Notice that the discrete-time
algorithm~\eqref{eq::DCDisc} is an equivalent iterative form of
\eqref{eq::DC1} obtained by Euler discretization with stepsize
$\delta$. When $\delta\to0$, we can expect that the stability and
convergence properties of~\eqref{eq::DCDisc} are similar to that of
\eqref{eq::DC1}, i.e., $x^i$ tracks the average of the network inputs
in its $O(\beta^{-1})$ neighborhood, provided the network topology is
strongly connected and weight-balanced digraph.  Notice that the
structure~\eqref{eq::DCDisc} allows us to circumvent discretizing the
derivative of the input signals and, as a result, avoid the one-step
delayed tracking reported in~\cite{MZ-SM:08a}. Next, note
that 
$u^i$ is never communicated directly.

Next, we explore the bounds on the stepsize $\delta$ such
that~\eqref{eq::DCDisc} is convergent and tracks the input
average. 
The proof of the results is presented in Appendix~\ref{app:A}.  We start by
studying the stability and convergence properties of the zero-system.

\begin{lem}[Convergence analysis and stepsize characterization of the
  zero-system of~\eqref{eq::DCDisc_smlr}]\label{lem::disc}
  Let $\GG$ be strongly connected and weight-balanced. For $\alpha$,
  $\beta>0$, the trajectory of the zero-system of discrete-time
  algorithm~\eqref{eq::DCDisc_smlr} over~$\GG$ starting from any
  initial condition $\vect{x}(0), \vect{v}(0) \in\real^N$ satisfies
  \begin{displaymath}x^i(k)\to
    -\frac{\alpha^{-1}}{N}\SUM{i=1}{N}v^i(0),\quad v^i(k)\to
    \avrg{v^j(0)}, \quad \forall i\in\until{N},
  \end{displaymath}
  asymptotically, as $k\to \infty$, provided
  $\delta\in(0,\min\{\alpha^{-1},\beta^{-1}(\text{d}_{\max}^{\text{out}})^{-1}\})$.\boxend
\end{lem} 

The following result establishes an upper bound on the solutions of
the algorithm~\eqref{eq::DCDisc} for any given initial conditions. In
the following, we let $\vect{\Phi}(k,j) =
(\vect{I}_{N-1}-\delta\beta\vect{R}^\top\vect{\mathsf{L}}\vect{R})^{k-j}$.

\begin{thm}[Upper bound on the tracking error
  of~\eqref{eq::DCDisc_smlr}]\label{thm::Disc_bound}
  Let $\GG$ be strongly connected and weight-balanced. Each agent has
  an input $u^i(k)$. For $\alpha$, $\beta>0$, the trajectory of the
  algorithm~\eqref{eq::DCDisc} over~$\GG$ starting from any initial
  condition $\vect{z}(0), \vect{v}(0) \in\real^N$ satisfies,
  \begin{align}\label{eq::xDCDisc_bound}
    &\left| x^i(k)-\avrg{u^j(k)}+\frac{1}{N}\delta
      \SUM{j=0}{k-1}(1-\delta\alpha)^j\SUM{j=1}{N}v^j(0)\right|\leq
    \Lnorm\vect{y}(t) + \delta\SUM{j=0}{k-1}(1-\delta\alpha)^j\vect{r}
    \vect{r}^\top\vect{w}(0)\Rnorm\leq\nonumber
    \\
    &\left|(1-\alpha\delta\SUM{j=0}{k-1}(1-\delta\alpha)^j)\right|
    \Lnorm \vect{y}(0)\Rnorm+\Lnorm\vect{\Phi}(k,0)
    \Rnorm\Lnorm\vect{y}(0)\Rnorm+\alpha
    \Lnorm(\SUM{j=0}{k-1}\vect{\Phi}(k-1,j)(1-\delta\alpha)^j)
    \Rnorm\Lnorm\vect{y}(0)\Rnorm+\nonumber
    \\
    &\Lnorm(\SUM{j=0}{k-1}\vect{\Phi}(k-1,j)(1-\delta\alpha)^j)
    \Rnorm\Lnorm\vect{w}(0)\Rnorm +
    \Lnorm\SUM{j=0}{k-1}\vect{\Phi}(k-1,j)(1-\delta\alpha)^j\Rnorm
    \Lnorm\Delta\vect{u}(0)\Rnorm+ \nonumber
    \\
    &\Lnorm\vect{R}\SUM{j=0}{k-1}\vect{\Phi}(k-1,j)
    \vect{R}^\top\Delta\vect{u}(j)\Rnorm,
  \end{align}
  for all $i \in \until{N}$, where $\vect{y}$ is defined
  in~\eqref{eq::xTOy} and $\vect{w}$ is
   \begin{equation}\label{eq::wTOv_disc}
     \vect{w}=\vect{v}-\Bvect{v},
     ~~~\Bvect{v}=\vect{\Pi}_N(\Delta
     \vect{u}(k)+\delta\alpha\vect{u}(k)).
  \end{equation}\boxend
\end{thm}

Next, we show that for networks with strongly connected and
weight-balanced digraph topologies, the discrete-time
algorithm~\eqref{eq::DCDisc} solves Problem~\ref{prob::main} with a
nonzero steady-state error, provided
$\delta\in(0,\min\{\alpha^{-1},\beta^{-1}(\text{d}_{\max}^{\text{out}})^{-1}\})$,
the algorithm is initialized properly and the essential norm of the
projection of the input difference vector into the agreement space is
bounded.

\begin{cor}[The algorithm~\eqref{eq::DCDisc} solves
  Problem~\ref{prob::main}]\label{cor::DCDiscProb1Sol}
  Let $\GG$ be strongly connected and weight-balanced.  Assume that
  the differences of the inputs of the network satisfy
  $\|\vect{\Pi}_N\Delta\vect{u}\|_{\text{ess}} = \gamma<\infty$. Then,
  for any $\alpha,\beta>0$, the algorithm~\eqref{eq::DCDisc}
  over~$\GG$ initialized at $z^i(0),v^i(0)\in \real$
  such that $\SUM{i=1}{N}v^i(0)=0$ solves Problem~\ref{prob::main} (in
  the output $x^i$) with an upper-bounded steady-state error provided $\delta \in(0,
  \min\{\alpha^{-1},\beta^{-1}(\text{d}_{\max}^{\text{out}})^{-1}\})$, specifically
  \begin{equation*}
    \lim_{k\to\infty} \left| x^i(k)-\avrg{u^j(k)}\right| \leq
    (\delta\beta\Hlambda_2)^{-1} \gamma,\quad  i \in \until{N}.
  \end{equation*}
\boxend
\end{cor}
One can make similar comments to those of
Remark~\ref{re:performance-regulation} regarding the tuning of the
performance of~\eqref{eq::DCDisc} via the design parameters $\alpha$
and $\beta$. In the following, we identify conditions, involving
inputs and their differences, under which the
algorithm~\eqref{eq::DCDisc} solves Problem~\ref{prob::main} with zero
steady-state error.

\begin{lem}[Conditions on inputs for zero steady-state error
  of~\eqref{eq::DCDisc}]\label{eq::Cnvrg_disc}
  Let $\GG$ be strongly connected and weight-balanced. Assume there
  exists $\delta\in(0,\min\{\alpha^{-1},
  \beta^{-1}(\text{d}_{\max}^{\text{out}})^{-1}\})$ and $\alpha>0$
  such that for all $i\in\until{N}$, one of the following conditions
  are satisfied
  \begin{itemize}
  \item[(a)] $\Delta u^i(k)+\delta\alpha u^i(k)$ 
    converges to a common dynamics $l(k)$;
  \item[(b)] $\Delta u^i(k+1)-\Delta u^i(k)+\delta\alpha \Delta
    u^i(k)$  converges to a common dynamics
    $l(k)$.
  \end{itemize}
  Then, the algorithm~\eqref{eq::DCDisc} over~$\GG$ with the given
  $\delta$ and $\alpha$, $z^i(0),v^i(0)\in\real$ such
  that $\SUM{i=1}{N}v^i(0)=0$, for any $\beta>0$, makes $x^i(k)\to
  \avrg{u^j(k)}$, as $k\to\infty$, for all $i \in \until{N}$.\boxend
\end{lem}

\vspace{-6pt}
\section{Dynamic average consensus with controllable rate of
  convergence and limited control
  authority}\label{sec::Rate}\vspace{-2pt}

In this section, we address the dynamic average consensus
Problems~\ref{prob::rate} and \ref{prob::sat}.  As discussed in
Section~\ref{sec::ProbDef}, the goal in setting up these problems is
to come up with an algorithm which is more suitable
for applications where the agreement state $x^i$
in~\eqref{eq::AgentSingInt} corresponds to some physical variable such
as position of a robotic system. In such networked systems, agents
might have limited control authority and can not implement the
high-rate commands dictated by the consensus algorithm.

Although the rate of convergence of the algorithm can be controlled by
the choice of $\alpha$ and $\beta$, these variables are centralized
variables and the effect is universal across the network. One can
expect that a more efficient consensus algorithm is one that allows
agents with limited power to move at their own pace. To this end, we
make a modification to the structure of the consensus
algorithm~\eqref{eq::DC1},
\begin{subequations}\label{eq::DC2}
  \begin{align}
    \dot{z}^i & =\dot{u}^i-\alpha
    (z^i-u^i)-\beta\SUM{j=1}{N}\vect{\mathsf{L}}_{ij}z^j-v^i,
    \label{eq::DC2-a}
    \\
    \dot{v}^i & =\alpha\beta\SUM{j=1}{N}\vect{\mathsf{L}}_{ij}z^j,
    \label{eq::DC2-aa}
    \\
    \dot{x}^i & =-\theta^i(t)(x^i-z^i) +\dot{u}^i-\alpha
    (z^i-u^i)-\beta\SUM{j=1}{N}\vect{\mathsf{L}}_{ij}z^j-v^i, \label{eq::DC2-b}
  \end{align}
\end{subequations}
where $\theta^i:[0,\infty) \to \real$ is a time-varying gain which is
bounded from below and above, i.e., at all $t\geq 0$ we have
$0<\underline{\theta}^i\leq \theta^i(t)\leq\bar{\theta}^i$, for $i \in
\until{N}$. As we show below, agents that wish to slow down their rate
of convergence use this gain to adjust it. Note the cascading
structure of the algorithm. As such, the stability properties
of~\eqref{eq::DC2-a}-\eqref{eq::DC2-aa} (\emph{information phase}) are
independent of~\eqref{eq::DC2-b} and are as characterized in
Section~\ref{sec::DyConsensus}. The information phase allows agents to
obtain the average with a convergence rate that is common across the
network. The dynamics~\eqref{eq::DC2-b} (\emph{motion phase}) allows
each agent $i \in \until{N}$ to tweak its convergence rate by
adjusting the gain $\theta^i$. 
We start our analysis by examining the rate of convergence of the
algorithm~\eqref{eq::DC2} and establishing an upper bound on its
tracking error.

\begin{lem}[The algorithm~\eqref{eq::DC2} solves
  Problem~\ref{prob::rate}]\label{lem::rate}
Let $\GG$ be strongly
  connected and weight-balanced. For inputs whose derivatives satisfy
  $\|\vect{\Pi}_N\dvect{u}\|_{\text{ess}} = \gamma<\infty$, for any
  $\alpha,\beta>0$ the algorithm~\eqref{eq::DC2} initialized
  at $x^i(0),v^i(0)\in \real$ such that
  $\SUM{i=1}{N}v^i(0)=0$, then we have the same ultimate tracking
  error bound of~\eqref{eq::DC1_ultimate_bound}. The rate of decay of
  the transient response is
  $\min\{\underline{\theta}^i,\alpha,\beta\Hlambda_2\}$ for each agent
  $i \in \until{N}$.
\end{lem}
\begin{pro}
  Consider the information
  phase~\eqref{eq::DC2-a}-\eqref{eq::DC2-aa}. From
  Theorem~\ref{thm::DC1_bound} and Corollary~\ref{cor::DC1Prob1Sol}, it
  follows that $z^i-\avrg{u^j(t)}$ has the ultimate bound
  \begin{equation}\label{eq::z_bound_CD2}
    \lim_{t\to\infty}\sup   \left| z^i(t)-\avrg{u^j(t)}\right|\leq
    (\beta\Hlambda_2)^{-1}\gamma,
  \end{equation}
  and converges to this neighborhood of the input average with a rate
  of $\min\{\alpha,\beta\Hlambda_2\}$.  Next, consider the motion
  phase~\eqref{eq::DC2-b}, which can be written as
  \begin{displaymath}
    \dot{x}^i=-\theta^i(t)(x^i-z^i) +\dot{z}^i,~~~~~ i \in
    \until{N},\forall t\geq0. 
  \end{displaymath}
  With the change of variables $d^i=x^i-z^i$, $i \in \until{N}$, this
  can be equivalently written as
  \begin{equation}\label{eq::DC2-b_equival}
    \dot{d}^i=-\theta^i(t)d^i,~~~~~ i \in \until{N},~\forall t\geq0.
  \end{equation}
  Using the Lyapunov function $V^i=\frac{1}{2}(d^i)^2$, it is not
  difficult to show that, for $0<\underline{\theta}^i\leq
  \theta^i(t)\leq\bar{\theta}^i$,~\eqref{eq::DC2-b_equival} is an
  exponentially stable system which satisfies the following bound
  \begin{displaymath}
    \left| x^i(t)-z^i(t)\right|=\left| d^i(t)\right|\leq\left| x^i(0)-z^i(0)
    \right|\e{-\underline{\theta}^i t},~~~~~ i \in \until{N},~\forall t\geq0.
  \end{displaymath}
  Therefore, 
  \begin{displaymath}
    \left| x^i(t)-\avrg{u^j(t)}\right|\leq\left| x^i(0)-z^i(0)
    \right|\e{-\underline{\theta}^i t}+\left|
      z^i(t)-\avrg{u^j(t)}\right|,\quad i \in \until{N},~\forall
    t\geq0. 
  \end{displaymath}
  Then, we conclude that~\eqref{eq::DC1_ultimate_bound} is
  satisfied. The rate of convergence of agent~$i$ is~$\min\{\underline{\theta}^i,\alpha,\beta\Hlambda_2\}$.
\end{pro}

As before, the design parameters $\alpha$ and $\beta$ can be used to
tune the overall rate of convergence.  Agents who wish to move at a
slower pace can use the motion phase with $\underline{\theta}^i\leq
\min\{\alpha,\beta\Hlambda_2\}$ to accomplish their goal.  The
time-varying nature of $\theta^i$ allows for agents to accelerate and
decelerate the convergence as desired. 
Notice that the ultimate error bound guaranteed by
algorithm~\eqref{eq::DC2} is the same as the one for
algorithm~\eqref{eq::DC1}. Therefore, the local first-order
filter~\eqref{eq::DC2-b} adjusts the rate of convergence without
having any adverse effect on the error bound.

\begin{rem}[Discrete-time implementation and switching networks]
  The results above can be extended to switching networks and
  discrete-time settings. For brevity this extension is omitted. In
  the discrete-time implementation, it is straightforward to show that
  for convergence we should require
  $\delta\in(0,\min\{\bar{\bar{\theta}}^{-1},\alpha^{-1},
  \beta^{-1}(\text{d}_{\max}^{\text{out}})^{-1}\})$, where
  $\bar{\bar{\theta}} = \underset{i\in \until{N}}{\max} \{
  \bar{\theta}^i\}$.\boxend
\end{rem}

Next, we consider the case when saturation is present in the driving
command.  The following result states that, under suitable conditions,
  the algorithm~\eqref{eq::DC2} is a solution for
Problem~\ref{prob::sat} with the same error bounds as if no saturation
was present.

\begin{lem}[The algorithm~\eqref{eq::DC2} solves
  Problem~\ref{prob::sat}]\label{lem::sat}
  Let $\GG$ be strongly connected and weight-balanced.  Suppose the
  driving command at each agent $i\in\until{N}$ is bounded by
  $\bar{c}^i>0$, i.e., $\dot{x}^i=\sat{\bar{c}^i}{c^i}$. Assume for
  every agent $i\in\until{N}$, the following holds: (a) the input
  signal at each agent is such that $\avrg{u^j}$ is bounded, the input
  derivatives satisfy $\|\vect{\Pi}_N\dvect{u}\|_{\text{ess}} =
  \gamma<\infty$, and $\|\dot{u}^i\|_{\text{ess}} = \mu^i<\infty$; (b)
  $\bar{c}^i>\mu^i+\gamma$. Then, for any $\alpha, \beta>0$, and
  constant $\theta^i>0$, the algorithm~\eqref{eq::DC2} starting from
  any $x^i(0),v^i(0)\in \real$ such that
  $\SUM{i=1}{N}v^i(0)=0$ satisfies that the ultimate tracking error
  bound~\eqref{eq::DC1_ultimate_bound}.
\end{lem}
\begin{pro}
  Following the proof
  of Lemma~\ref{lem::rate}, for the information phase~\eqref{eq::DC2-a} and~\eqref{eq::DC2-aa},  we
  have~\eqref{eq::z_bound_CD2}. To complete the proof, we will show
  that under the given conditions for the input signals, despite the
  saturation, $x^i\to z^i$ asymptotically for all
  $i\in\until{N}$. Under the saturation constraint,~\eqref{eq::DC2-b}
  takes the form $\dot{x}^i =-\sat{\bar{c}^i}{\theta^i(x^i-z^i)
    +\dot{z}}$, for $i\in\until{N}$.  The rest of the proof relays on
  Proposition~\ref{pro::sat2}. According to this result, we need to
  show that a) $z^i$ is a bounded signal; b)
  $|\dot{z}^i(t)|<\bar{c}^i$ for all $t>t^\star$ where $t^\star$ is
  some finite time. For any given finite initial conditions and input
  signals with bounded average the requirement (a) is satisfied due to
  convergence guarantees of \eqref{eq::DC2-a}-\eqref{eq::DC2-aa}. In
  the following, we show that the requirement (b) is also satisfied
  due to the given assumptions.  With change of
  variables~\eqref{eq::wTOv} and
  $\vect{y}=\vect{z}-\avrg{u^j}\vect{1}_N$, we can
  represent~\eqref{eq::DC2-a} as $\dvect{z}
    =-\alpha\vect{y}-\beta\vect{\mathsf{L}}\vect{y}-\vect{w} +
    \avrg{\dot{u}^j}\vect{1}_N$.~Therefore,
  \begin{displaymath}
    \lim_ {t\to\infty}\left|
      \dot{z}^i(t)\right|\leq   \lim_ {t\to\infty}\left|
    -\alpha\vect{y}^i(t)-\vect{w}^i(t) +
    \avrg{\dot{u}^j(t)}\right| + \lim_ {t\to\infty}\Lnorm
    \beta\vect{\mathsf{L}}\vect{y}(t)\Rnorm \quad i\in\until{N}.
  \end{displaymath}
Using the results and the variables introduced in the proof of
  Theorem~\ref{thm::DC1_bound}, we can show that
  \begin{displaymath}
    -\alpha\vect{y}-\vect{w}+\avrg{\dot{u}^j}\vect{1}_N = 
    -\begin{bmatrix} \alpha\vect{S}_{11}+\vect{S}_{21}
      &\alpha\vect{S}_{12}+\vect{S}_{22} 
    \end{bmatrix}\begin{bmatrix}
      \vect{y}(0)\\
      \vect{w}(0)
    \end{bmatrix}-\e{-\alpha
      t}\vect{R}\vect{R}^\top\dvect{u}(0)+\dvect{u}(t),
  \end{displaymath} 
  where $\vect{S}_{11}$ and $\vect{S}_{12}$ are given in~\eqref{eq::S},
  and we have
  \begin{align}\label{eq::S2p}
    &\vect{S}_{21} = -\alpha
    \vect{R}\vect{\Phi}(t,0)+\alpha^2\vect{R}(\int_{0}^t
    \vect{\Phi}(t,\tau)\e{-\alpha\tau}d\tau)\vect{R}^\top +
    \alpha\vect{R}\vect{R}^\top\e{-\alpha t},\nonumber
    \\
    &\vect{S}_{22} = \vect{r}\vect{r}^\top+ \alpha
    \vect{R}(\int_{0}^t\vect{\Phi}(t,\tau)\e{-\alpha\tau}d\tau)\vect{R}^\top
    + \vect{R}\vect{R}^\top\e{-\alpha t}.
  \end{align}
 Recall that $\vect{\Phi}(t,\tau) =
  \e{-\beta\vect{R}^\top\vect{\mathsf{L}}\vect{R}(t-\tau)}$, then,
  \begin{align*}
    & \Lnorm \beta\vect{\mathsf{L}} \vect{R}\int_{0}^t\vect{\Phi}(t,\tau)
    \vect{R}^\top\dvect{u}(\tau)d\tau\Rnorm =
    \Lnorm\beta\vect{R}\vect{R}^\top\vect{\mathsf{L}}\vect{R}\int_{0}^t\vect{\Phi}(t,\tau)
    \vect{R}^\top\dvect{u}(\tau)d\tau\Rnorm =
    \\
    &\Lnorm \vect{R} \e{-\beta\vect{R}^\top\vect{\mathsf{L}}\vect{R}t}
    \int_{0}^t \beta\vect{R}^\top\vect{\mathsf{L}}\vect{R}
    \e{\beta\vect{R}^\top\vect{\mathsf{L}}\vect{R}\tau} \vect{R}^\top\vect{\Pi}_N
    \dvect{u}(\tau)d\tau\Rnorm\leq
    \\
    &\Lnorm\vect{R} \e{-\beta\vect{R}^\top\vect{\mathsf{L}}\vect{R}t}
    \int_{0}^t \beta\vect{R}^\top\vect{\mathsf{L}}\vect{R}
    \e{\beta\vect{R}^\top\vect{\mathsf{L}}\vect{R}\tau} d\tau
    \|\vect{\Pi}_N\dvect{u}\|_{\text{ess}}\Rnorm=
    \\
    &\Lnorm\vect{R} \e{-\beta\vect{R}^\top\vect{\mathsf{L}}
      \vect{R}t}(\e{\beta\vect{R}^\top\vect{\mathsf{L}}\vect{R}t}-\vect{I}_N)
    \|\vect{\Pi}_N\dvect{u}\|_{\text{ess}}
    \Rnorm\leq\|\vect{\Pi}_N\dvect{u}\|_{\text{ess}}+\e{-\beta\Hlambda_2 t}
    \|\vect{\Pi}_N\dvect{u}\|_{\text{ess}}
  \end{align*}
  Recall~\eqref{eq::SolDC1_y}. In light of the relations above 
  we can show that
  \begin{displaymath}
    \lim_ {t\to\infty}\left|
      \dot{z}^i(t)\right| \leq \mu^i+\gamma,\quad i\in\until{N}.
  \end{displaymath}
  Therefore, there exists a finite time $t^\star$ such that
  $|\dot{z}^i(t)|<\bar{c}^i$ for all $t>t^\star$ and $i\in\until{N}$.
\end{pro}

\vspace{-6pt}
\section{Dynamic average consensus with privacy
  preservation}\label{sec::Priv}\vspace{-2pt}

Here, we study the dynamic average consensus problem with privacy
preservation. We consider adversaries that do not interfere with the
implementation of the algorithm but are interested in retrieving
information about the inputs, their average, or the agreement state
trajectories of the individual agents. These adversaries might be
\emph{internal}, i.e., part of the network, or
\emph{external}. Internal adversaries have access at no cost to
certain information that external adversaries do not. More
specifically, an internal adversary has knowledge of the parameters
$\alpha$, $\beta$ of the algorithm~\eqref{eq::DC1}, its corresponding
row in the Laplacian matrix, and the agreement state of its
out-neighbors. We also assume that the agent is aware of whether the
algorithm is initialized with $\vect{v}(0)=\vect{0}$.  We refer to the
extreme case when an internal adversary knows the whole Laplacian
matrix and the initial conditions of its out-neighbors as a
\emph{privileged internal adversary}.  Regarding external adversaries,
we assume they have access to the time history of all the
communication messages.  We refer to the extreme case when an external
adversary has additionally knowledge of the parameters $\alpha$,
$\beta$, the Laplacian matrix, and the initial conditions as a
\emph{privileged external adversary}.

The next result characterizes the privacy-preservation properties of
the dynamic average consensus algorithm~\eqref{eq::DC1} against
adversaries. Specifically, we show that this algorithm satisfies
Problem~\ref{prob::priv}(a).

\begin{lem}[The algorithm~\eqref{eq::DC1} preserves the privacy of the
  local inputs against adversaries]\label{lem::privLocalu}
  Let $\GG$ be strongly connected and weight-balanced.  The executions
  of the algorithm~\eqref{eq::DC1} over $\GG$ with $\alpha$,
  $\beta>0$, initialized at $x^i(0),v^i(0)\in\real$ such
  that $\sum_{i=1}^{N}v^i(0)=0$, satisfy
  \begin{itemize}
  \item[(a)] an external (respectively internal) adversary cannot
    reconstruct the input of any (respectively another) agent;
  \item[(b)] a privileged adversary cannot reconstruct the input of
    agent $i \in \until{N}$ as long as there exists $\bar{t}>0$ such
    that $\dot{u}^i(t)\neq 0$ for $t\in[0,\bar{t})$.
  \end{itemize}
\end{lem}
\begin{pro}
  First, we investigate the validity of claim (a). Using the results
  in the proof of Theorem~\ref{thm::DC1_bound} and recalling the
  change of variables~\eqref{eq::xvToyw}, the solution of the
  algorithm~\eqref{eq::DC1} for given initial conditions
  $x^i(0),v^i(0)\in\real$, for $i\in\until{N}$ can be
  written as follows
  \begin{align}\label{eq::FullSolCD1} 
    &\begin{bmatrix}
      \vect{x}(t)\\
      \vect{v}(t)
    \end{bmatrix}=
    \begin{bmatrix}
      \vect{S}_{11}&\vect{S}_{12}
      \\
      \vect{S}_{21}&\vect{S}_{22}
    \end{bmatrix}
    \begin{bmatrix}
      \vect{x}(0)+(\avrg{u^j(0}))\vect{1}_N
      \\
      \vect{v}(0)+\vect{\Pi}_N(\dvect{u}(0)+\alpha\vect{u}(0))
    \end{bmatrix}+\begin{bmatrix}
      (\avrg{u^j(t)})\vect{1}_N
      \\
      \vect{\Pi}_N(\dvect{u}(t)+\alpha\vect{u}(t))
    \end{bmatrix}
    +
    \\
    &\begin{bmatrix} -\vect{R}\int_{0}^t\vect{\Phi}(t,\tau)
      \e{-\alpha\tau}d\tau\,\vect{R}^\top\dvect{u}(0) +
      \vect{R}\int_{0}^t\vect{\Phi}(t,\tau)
      \vect{R}^\top\dvect{u}(\tau)d\tau
      \\
      \alpha\vect{R}\int_{0}^t\vect{\Phi}(t,\tau)
      \e{-\alpha\tau}d\tau\,\vect{R}^\top\dvect{u}(0)-\alpha\vect{R}\int_{0}^t
      \vect{\Phi}(t,\tau) \vect{R}^\top\dvect{u}(\tau)d\tau+\e{-\alpha
        t}\vect{R}\vect{R}^\top\dvect{u}(0)-\vect{R}\vect{R}^\top\dvect{u}(t)
    \end{bmatrix},
    \notag
  \end{align}
  where $\vect{S}_{11}$ and $\vect{S}_{12}$ are given
  in~\eqref{eq::S}, and $\vect{S}_{21}$ and $\vect{S}_{22}$ are given
  in~\eqref{eq::S2p}.  For an external adversary that only has
  knowledge of the time history of $\vect{x}$, the number of unknowns
  in~\eqref{eq::FullSolCD1} (i.e., $\vect{u}(0)$, $\vect{u}(t)$,
  $\dvect{u}(t)$, $\vect{v}(t)$, for $\forall t\geq0$, $\alpha$,
  $\beta$ and $\vect{\mathsf{L}}$), regardless of the initial condition
  requirement $\sum_{i=1}^{N}v^i(0)=0$, is larger than the number of
  equations. This is true even if the inputs are static. Thus, the
  claim (a) for external adversaries follows.  Regarding the claim (a)
  for internal adversaries, we consider the extreme case where the
  adversarial agent, say $j$, is the in-neighbor of every other agent
  in the network, and therefore knows the time history of the
  aggregated vector $\vect{x}$. Now consider~\eqref{eq::DC1-b} for all
  $i\in\VV \setminus \{j\}$. Recall that agent $j$ does not know
  $\vect{\mathsf{L}}_{ik}$, $k\in\VV$, of all agent $i\in\VV \setminus
  \{j\}$. Therefore, even if it knows the initial condition $v^i(0)$,
  it cannot obtain $v^i(t)$, $t>0$. Next consider~\eqref{eq::DC1-a},
  and again assume an extreme case that the adversarial agent $j$ can
  numerically reconstruct $\dot{x}^i$ with an acceptable precision and
  the inputs are static. Despite these assumptions, because $u^i$ and
  $\sum_{k=1}^{N}\vect{\mathsf{L}}_{ik}x^k$, $\forall t\geq 0$ of all agent
  $i\in\VV \setminus \{j\}$ are unknown to agent $j$, regardless of
  value of $v^i$, this agent cannot reconstruct $u^i$
  from~\eqref{eq::DC1-b}. This concludes validity of the claim (a) for
  internal agents.
  
  Next, we examine claim (b) considering both the internal and
  external adversary case at the same time.  For an internal
  adversary, assume the extreme case when it is the in-neighbor of
  every other agent in the network. As a result, it knows the time
  history of the aggregated vector $\vect{x}$. At any given $\tau>0$,
  using its knowledge of $\vect{x}(t)$ over $t\in[0,\tau]$ and the
  information on the initial conditions and the parameters of the
  algorithm, a privileged internal or external adversary can
  reconstruct $v^i(t)$, $i\in\until{N}$, for all $t\in[0,\tau]$ by
  integrating \eqref{eq::DC1-b}. The adversary can also use its
  knowledge of $\vect{x}(t)$ over $t\in[0,\tau]$ to construct
  numerically $\dvect{x}(t)$ over the same period of time. Then, the
  adversary using~\eqref{eq::DC1-a}, knows the right-hand side of the
  following equation
  \begin{equation}\label{eq::DC1-a_uODE}
    \dot{u}^i+\alpha u^i=-\dot{x}^i -\alpha
    x^i-\beta\SUM{j=1}{N}\vect{\mathsf{L}}_{ij}x^j-v^i,\quad\forall i\in\until{N}.
  \end{equation}
  Because there exists $\bar{t}>0$ such that $\dot{u}^i(t)\neq 0$ for
  $t\in[0,\bar{t})$, \eqref{eq::DC1-a_uODE} is an ordinary
  differential equation (ODE) with variable $u^i$. The adversary does
  not know the initial condition $u^i(0)$, hence, it cannot obtain the
  unique solution of the ODE, i.e., the dynamic input $u^i$. This
  validates claim (b).
\end{pro}

\begin{rem}[Privacy preservation of static inputs against privileged
  adversaries]\label{rem::static_priv_DC1}
  To protect local static inputs from privileged adversaries, agents
  can add a static or time-varying value to their inputs at the
  beginning for some short period of time (so that the requirement of
  Lemma~\ref{lem::privLocalu}(b) is satisfied) and then remove
  it. This modification does not affect the final convergence
  properties of the algorithm~\eqref{eq::DC1}. 
  \boxend
\end{rem}

In general, the algorithm~\eqref{eq::DC1} does not satisfy the
requirements (b) and (c) of Problem \ref{prob::priv}. Here, we propose
a slight extension of~\eqref{eq::DC2} that overcomes this
shortcoming. For each $i\in\until{N}$,~let
\begin{subequations}\label{eq::DC3}
  \begin{align}
    &\dot{z}^i =\dot{u}^i -\alpha (z^i-u^i) -
    \beta\SUM{j=1}{N}\vect{\mathsf{L}}_{ij}\tilde{z}^j-v^i ,\label{eq::DC3-a}
    \\
    &\dot{v}^i =
    \alpha\beta\SUM{j=1}{N}\vect{\mathsf{L}}_{ij}\tilde{z}^j,\label{eq::DC3-aa}
    \\
    &\dot{x}^i = -\theta^i(t)(x^i-z^i)+\dot{u}^i-\alpha(z^i-u^i) -
    \beta\SUM{j=1}{N}\vect{\mathsf{L}}_{ij}\tilde{z}^j-v^i ,\label{eq::DC3-b}
    \\
    &\tilde{z}^i=z^i+\psi(t),\label{eq::DC3-aaa}
  \end{align}
\end{subequations}
where $\psi:[0,\infty)\to\real$ is a common dynamic signal which is
known to all agents. Also, $\theta^i: [0,\infty)\to\real$ such that
$\underline{\theta}^i\leq \theta^i(t)\leq \bar{\theta}^i$ for all
$t\geq 0$ is a local signal only known to agent $i$.  The role of the
signal $\psi$ is to conceal the final agreement value from the
external adversaries to satisfy the item (b) in
Problem~\ref{prob::priv}.  Note that, because
$\SUM{j=1}{N}\vect{\mathsf{L}}_{ij}=0$, the signal $\psi$ has no effect on the
algorithm execution, and therefore, the executions of
algorithms~\eqref{eq::DC3} and~\eqref{eq::DC2} are the
same. Consequently, Lemma~\ref{lem::rate} is valid for~\eqref{eq::DC3}
as well.  As agents communicate $\tilde{z}^i$ instead of $z^i$, and
the signal $\psi$ is unknown to the external adversaries, recovering
the steady-state solution of the algorithm is impossible for such
adversaries.  The agreement state equation of any agent $i$
in~\eqref{eq::DC3-b} is a local equation, with all the components set
by that agent.  Therefore, $x^i(0)$ and $\theta^i$ can easily be
concealed from other agents, making it impossible for adversaries to
reconstruct the trajectories of $x^i$. This allows us to satisfy the
item (c) in Problem~\ref{prob::priv}.  The following result shows that
the algorithm~\eqref{eq::DC3} is privacy preserving and solves
Problem~\ref{prob::priv}.  Its proof is a consequence of the above
discussion and Lemmas~\ref{lem::rate} and~\ref{lem::privLocalu}, and
is omitted for brevity.
\begin{lem}[The algorithm~\eqref{eq::DC3} solves
  Problem~\ref{prob::priv}]\label{lem::solutionProbPriv} Under the
  hypotheses of Lemma~\ref{lem::rate}, the ultimate tracking error
  bound~\eqref{eq::DC1_ultimate_bound} is valid for all trajectories
  $t \mapsto x^i(t)$ of the algorithm~\eqref{eq::DC3}. Furthermore,
  \begin{itemize}
  \item[(a)] an external (respectively internal) adversary cannot
    reconstruct the input of any (respectively another) agent;
  \item[(b)] a privileged adversary cannot reconstruct the input of
    agent $i \in \until{N}$ as long as there exists $\bar{t}>0$ such
    that $\dot{u}^i(t)\neq 0$ for $t\in[0,\bar{t})$;
  \item[(c)] external adversaries cannot obtain the final agreement
    value of the network as long as $\psi$ is unknown to them;
  \item[(d)] an adversary cannot reconstruct the trajectory $t \mapsto
    x^i(t)$ of agent $i \in \until{N}$ as long as $x^i(0)$ or
    $\theta^i$ is unknown to it.
  \end{itemize}
\end{lem}

\vspace{-6pt}
\section{Simulations}\label{sec:simulations}\vspace{-2pt}

\begin{figure}[t]
\centering
  \captionsetup[subfloat]{captionskip=10pt}
    \subfloat [] 
{%
      \begin{tikzpicture}[auto,thick,scale=1, every node/.style={scale=1}]
          \tikzstyle{mynode}=%
          [%
            minimum size=12pt,%
            inner sep=0pt,%
            outer sep=0pt,%
            ball color=red!20!orange!70,
            shape=circle%
          ]
          
  \draw
(-0.5,0) node[mynode] (1) {{\scriptsize1}}
    (0.5,0) node[mynode] (2) {{\scriptsize2}}
    (1,-1) node[mynode] (3) {{\scriptsize3}}
     (0.5,-2) node[mynode] (4) {{\scriptsize4}}
     (-0.5,-2) node[mynode] (5) {{\scriptsize5}}
     (-1,-1) node[mynode] (6) {{\scriptsize6}};
           \draw[-latex]  (1)->(2) ;
      \draw[-latex]  (2)->(3) ;
     \draw[-latex]  (3)->(4) ;
     \draw[-latex] (4)->(5) ;
          \draw[-latex] (5)->(6) ;
               \draw[-latex]  (6)->(1) ;
\end{tikzpicture}
    }\quad
\subfloat [] {%
  \begin{tikzpicture}[auto,thick,scale=1, every node/.style={scale=1}]
          \tikzstyle{mynode}=%
          [%
            minimum size=12pt,%
            inner sep=0pt,%
            outer sep=0pt,%
            ball color=red!20!orange!70,%
            shape=circle%
          ]
          \draw
  (-0.5,0) node[mynode] (1) {{\scriptsize1}}
    (0.5,0) node[mynode] (2) {{\scriptsize2}}
    (1,-1) node[mynode] (3) {{\scriptsize3}}
     (0.5,-2) node[mynode] (4) {{\scriptsize4}}
     (-0.5,-2) node[mynode] (5) {{\scriptsize5}}
     (-1,-1) node[mynode] (6) {{\scriptsize6}};
           \draw[-latex]  (1)->(2) ;
           \draw[-latex]  (2)->(6) ;
           \draw[-latex]  (6)->(1) ;
     \draw[-latex]  (3)->(5) ;
     \draw[-latex] (5)->(4) ;
          \draw[-latex] (4)->(3) ;
          \end{tikzpicture}
    } \quad
\subfloat [] {%
   \begin{tikzpicture}[auto,thick,scale=1, every node/.style={scale=1}]
          \tikzstyle{mynode}=%
          [%
            minimum size=12pt,%
            inner sep=0pt,%
            outer sep=0pt,%
            ball color=red!20!orange!70,%
            shape=circle%
          ]
  \draw        
  (-0.5,0) node[mynode] (1) {{\scriptsize1}}
    (0.5,0) node[mynode] (2) {{\scriptsize2}}
    (1,-1) node[mynode] (3) {{\scriptsize3}}
     (0.5,-2) node[mynode] (4) {{\scriptsize4}}
     (-0.5,-2) node[mynode] (5) {{\scriptsize5}}
     (-1,-1) node[mynode] (6) {{\scriptsize6}};
           \draw[latex'-latex'] (3) --  (2);   
\end{tikzpicture}
   }  \quad
\subfloat [] 
{%
      \begin{tikzpicture}[auto,thick,scale=1, every node/.style={scale=1}]
          \tikzstyle{mynode}=%
          [%
            minimum size=12pt,%
            inner sep=0pt,%
            outer sep=0pt,%
            ball color=red!20!orange!70,%
            shape=circle%
          ]
          
  \draw
(-0.5,0) node[mynode] (1) {{\scriptsize1}}
    (0.5,0) node[mynode] (2) {{\scriptsize2}}
    (1,-1) node[mynode] (3) {{\scriptsize3}}
     (0.5,-2) node[mynode] (4) {{\scriptsize4}}
     (-0.5,-2) node[mynode] (5) {{\scriptsize5}}
     (-1,-1) node[mynode] (6) {{\scriptsize6}};
            \draw[-latex]  (1)->(2) ;
           \draw[-latex]  (2)->(6) ;
           \draw[-latex]  (6)->(1) ;
           \end{tikzpicture}
    }\quad
\subfloat [] 
{%
      \begin{tikzpicture}[auto,thick,scale=1, every node/.style={scale=1}]
          \tikzstyle{mynode}=%
          [%
            minimum size=12pt,%
            inner sep=0pt,%
            outer sep=0pt,%
            ball color=red!20!orange!70,%
            shape=circle%
          ]
          
  \draw
(-0.5,0) node[mynode] (1) {{\scriptsize1}}
    (0.5,0) node[mynode] (2) {{\scriptsize2}}
    (1,-1) node[mynode] (3) {{\scriptsize3}}
     (0.5,-2) node[mynode] (4) {{\scriptsize4}}
     (-0.5,-2) node[mynode] (5) {{\scriptsize5}}
     (-1,-1) node[mynode] (6) {{\scriptsize6}};
           \draw[-latex]  (3)->(4) ;
     \draw[-latex] (4)->(5) ;
          \draw[-latex] (5)->(3) ;
                     \draw[latex'-latex'] (5) --  (6);   
          
\end{tikzpicture}
    }\\
    \caption{Weight-balanced digraphs used in simulation (all edge
      weights are equal to $1$).}\label{fig::network}
\end{figure}
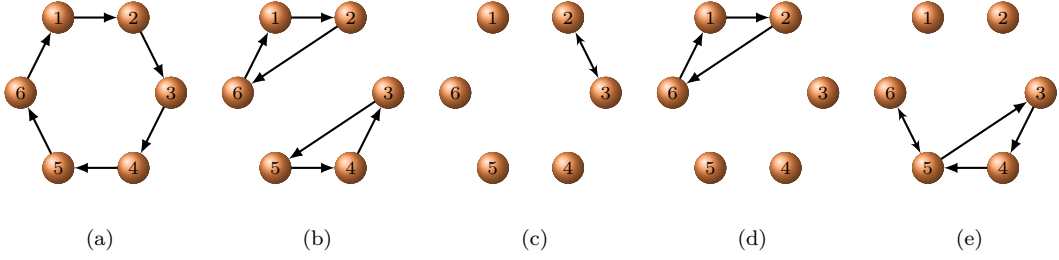

Here, we evaluate the performance of the proposed dynamic average
consensus algorithms in a number of scenarios. Fig.~\ref{fig::network}
shows the weight-balanced digraphs employed in the simulation.

\vspace{-6pt}
\subsection{Networks with time-varying interaction
  topologies}\label{sec::NExam_Swtch}\vspace{-2pt}

\begin{figure}[t!]
  \unitlength=0.5in \centering
  \subfloat[Case 1]{
    \includegraphics[width=0.45\textwidth]{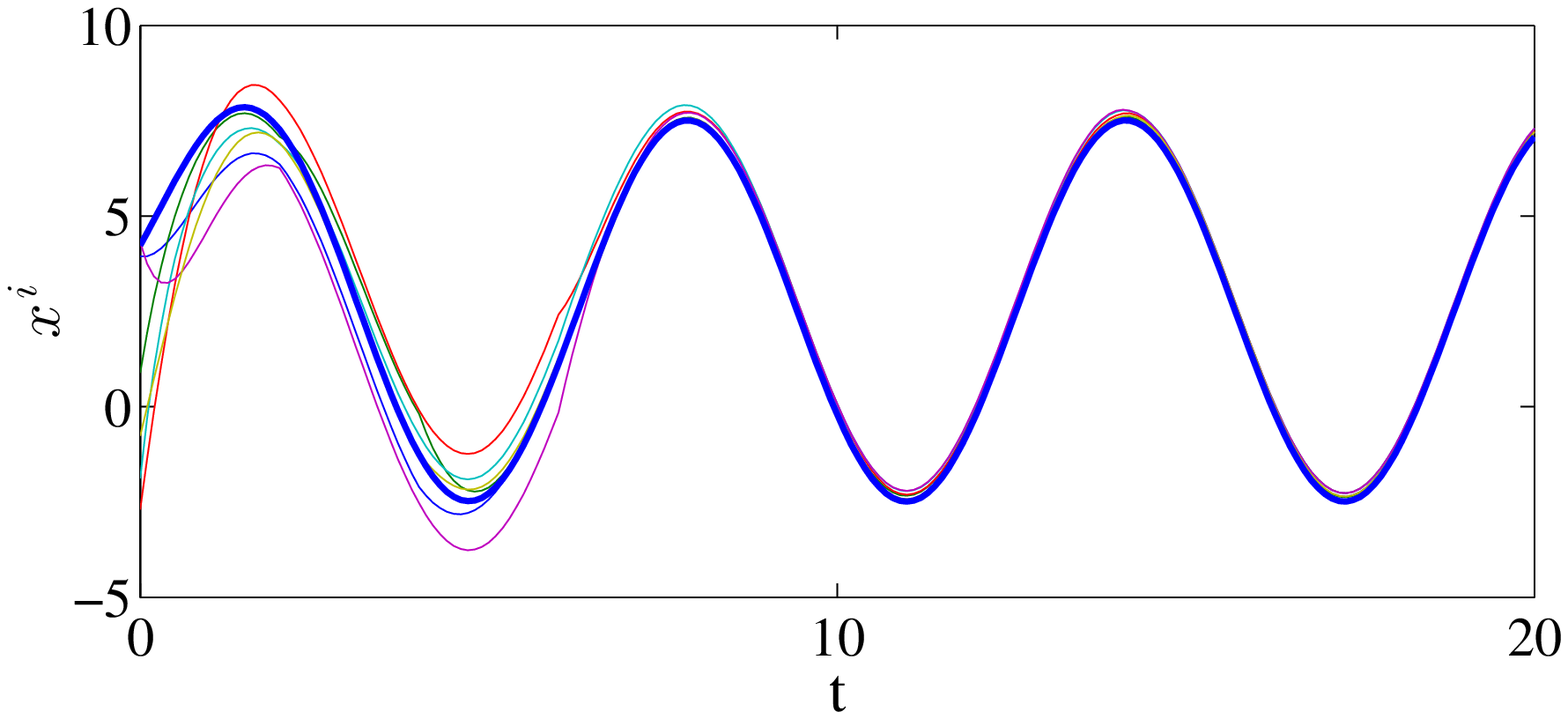}
  }\quad
  \subfloat[Case 2]{
    \includegraphics[width=0.45\textwidth]{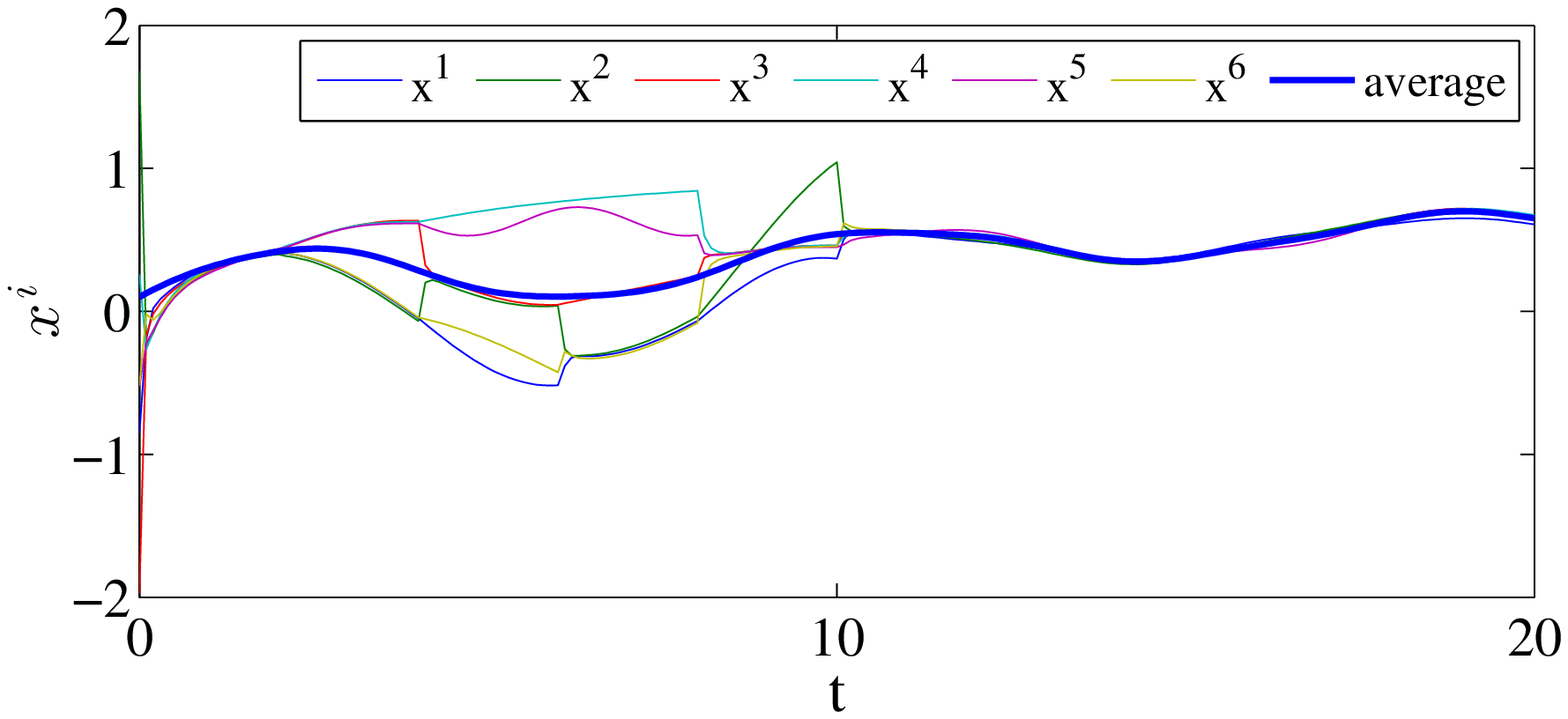} 
  }
  \caption{Simulation results for Case 1 and Case 2 of the numerical
    example of Section~\ref{sec::NExam_Swtch}: Solid thick blue line
    (colored thin lines)
    is the input average (resp. agreement state of
    agents). 
  }\label{fig::Ex2sim_Dynam}
\end{figure}

Consider a group of $6$ agents whose communication topology is
time-varying. We consider the following cases for the input signals

\begin{displaymath}
  \text{Case 1:~}
\left\{\begin{array}{l}
    u^1(t)=5\sin t+\frac{1}{t+2}+3,\\
    u^2(t)=5\sin t+\frac{1}{(t+2)^2}+4,\\
    u^3(t)=5\sin t+\frac{1}{(t+2)^3}+5,\\
    u^4(t)=5\sin t+10\e{-t}+4,\\
    u^5(t)=5\sin t+\mathrm{atan}\, t-1.5,\\
    u^6(t)=5\sin t-\tanh t+1.
  \end{array}\right.\quad\quad\text{Case 2:~}\left\{\begin{array}{l}
    u^1(t)=0.55\sin(0.8t),\\
    u^2(t)=0.5\sin (0.7t)+0.5\cos(0.6 t),\\
    u^3(t)=0.1 t,\\
    u^4(t)=\mathrm{atan}(0.5t),\\
    u^5(t)=0.1\cos(2t),\\
    u^6(t)=0.5\sin (0.5t).
  \end{array}\right.
\end{displaymath}
In Case 1, the communication topology iteratively changes, in
alphabetical order, every two seconds among the digraphs in
Fig.~\ref{fig::network}(b)-(e). In Case 2, the communication topology
changes, in alphabetical order, every two seconds among the digraphs
in Fig.~\ref{fig::network}(a)-(e). After $t=10$ seconds, the
communication topology is fixed at the digraph in
Fig.~\ref{fig::network}(a). Figure~\ref{fig::Ex2sim_Dynam} shows the
simulation results generated by implementing the
algorithm~\eqref{eq::DC1} with the following parameters: in Case~1,
$\alpha=\beta=1$ and in Case~2, $\alpha=3$ and $\beta=10$.

These examples show that, as long as the switching signal belongs to
$\mathcal{S}_{\mathrm{admis}}$, the agreement state $x^i$ stays
bounded. In Case 1, because the input signals converge to a common
function, the version of Lemma~\ref{eq::Cnvrg} for switching networks
implies that the algorithm~\eqref{eq::DC1} converges to the average
with zero steady-state error. However, in Case 2, we only can
guarantee tracking with bounded steady-state error. During the times
that the network is only weight-balanced, the error grows but still
stays bounded. One can expect that each connected group converges to
their respective input average. During these periods of time, there is
no way for separate components to have knowledge of the other groups'
inputs. However, once the network is strongly connected and
weight-balanced, then~\eqref{eq::DC1} resumes its tracking of the
input average across all network, as expected.

\vspace{-6pt}
\subsection{Dynamic inputs offset by a static
  value}\label{sec::sim_disc}\vspace{-2pt}

\begin{figure}[t!]
\unitlength=0.5in
\centering
 \includegraphics[width=0.75\textwidth]{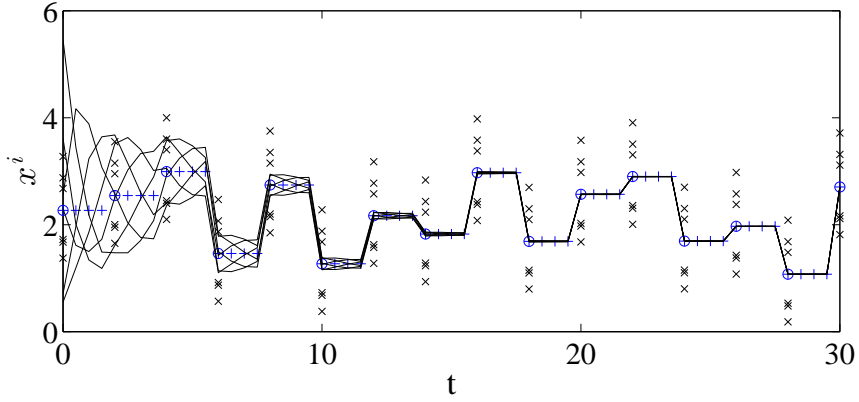}
 \caption{Simulation results for the numerical example of
    Section~\ref{sec::sim_disc}; The solid lines: the agreement states
   of~\eqref{eq::DCDisc};  $\times$: sampling points at $m\Delta
   t$; $\circ$: the average at $m\Delta t$; $+$: the average at
   $k\delta$. }\label{fig::Ex1sim_disc}
\end{figure}

Consider a process described by a fixed value plus a sine wave whose
frequency and phase are changing randomly over time. A group of $6$
agents with the communication topology shown in
Fig.~\ref{fig::network}(a) monitors this process by taking synchronous
samples, each according to
\begin{align*}
  u^i(m) = 2 + \sin(\omega(m) t(m)+\phi(m)) + b^i, ~~~m=0,1,\dots .
\end{align*}
Because of the unknown fixed bias $b^i$ of each agent, after each
sampling, every agent wants to obtain the average of the measurements
across the network before the next sampling time. Here, $\omega\sim
\mathrm{N}(0,0.25)$, $\phi\sim \mathrm{N}(0,(\pi/2)^2)$, with
$\mathrm{N}(.,.)$ indicating a Gaussian distribution. The data is
sampled at $0.5$ Hertz, i.e., $\Delta t=2$ seconds. The bias at each
agent is $b^1= -0.55$, $b^2=1$, $b^3=0.6$, $b^4=-0.9$, $b^5=-0.6$, and
$b^6=0.4$. Between sampling times $m$ and $m+1$, the input $u^i(k)$ is
fixed at $u^i(m)$.  Figure~\ref{fig::Ex1sim_disc} shows the result of
the simulation using the discrete-time~consensus~algorithm~\eqref{eq::DCDisc} with $\alpha=\beta=1$. The communication
bandwidth is 2 Hertz, i.e., $\delta=0.5$ seconds.  The application
of~\eqref{eq::DCDisc} results in perfect tracking after some time as
forecasted by Lemma~\ref{eq::Cnvrg_disc}.  Notice that here as it is impossible for the agents to know $u^i(-1)$, the use of the
algorithm in~\cite{MZ-SM:08a}, which requires the agents to initialize
their agreement states at $u^i(-1)$, results in tracking with a
steady-state error.

\vspace{-6pt}
\subsection{Limited control
  authority}\label{sec::NExam_Sat}\vspace{-2pt}

\begin{figure}[t!]
  \unitlength=0.5in \centering
  \subfloat[Dynamic average consensus algorithm~\eqref{eq::DC1} ] {
    \includegraphics[width=.45\linewidth]{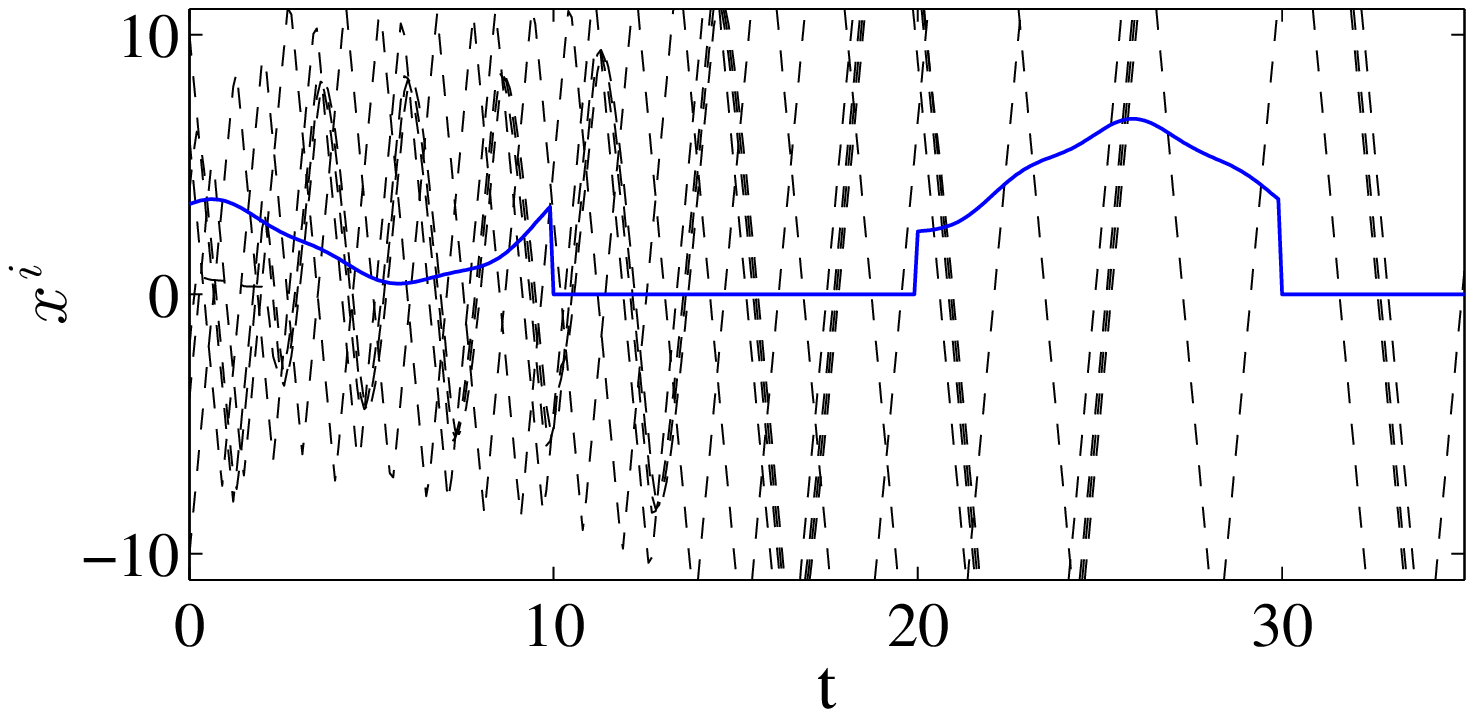}
  } \quad
  \subfloat[Dynamic average consensus algorithm~\eqref{eq::DC2}] {
    \includegraphics[width=.45\textwidth]{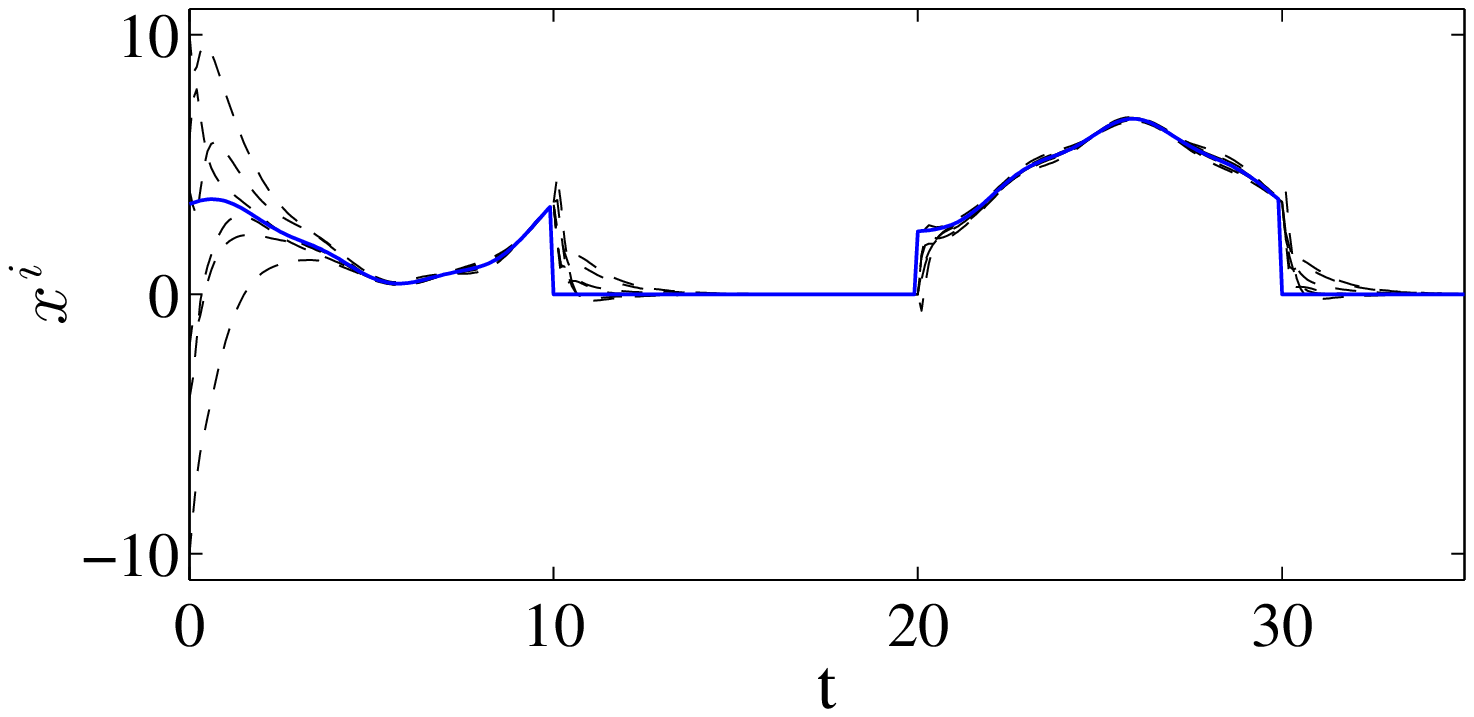} 
  }
  \caption{Simulation results for the numerical example of
    Section~\ref{sec::NExam_Sat}: Solid blue line (black dashed lines)
    is the input average (resp. agreement state of
    agents).}\label{fig::Ex3sim_Sat}
\end{figure}

We use the following numerical example to demonstrate the performance
of the algorithms~\eqref{eq::DC1} and~\eqref{eq::DC2} when the
driving command is bounded. Consider a group of 6 agents whose
communication topology is given in Fig.~\ref{fig::network}(a). The
input signals are as follows
\begin{displaymath}
  \begin{array}{ll}u^1(t)=u(t)\,(4\,\cos(0.5t)+10),&
    u^2(t)=u(t)(4\,\text{tanh}(t-5)+4\,\text{tanh}(t-25)+5),
    \\
    u^3(t)=\,u(t)(4\,\sin(0.5t+1)+8),&u^4(t)=u(t)(4\,\text{atan}(0.5t-5)-6),
    \\
    u^5(t)=\,u(t)(\sin(2t)-5),&
    u^6(t)=u(t)(4\,\cos(0.5t)+7),\end{array}
\end{displaymath}
where $u(t)=\sum_{i=0}^\infty((-1)^iH(t-10\,i))$, in which $H$ is the
step function, $H(t)=0$ if $t<0$, and $H(t) = 1$ if $t\geq 0$. For
both algorithms~\eqref{eq::DC1} and~\eqref{eq::DC2} we use $\alpha=10$
and $\beta=15$. In the algorithm~\eqref{eq::DC2} we set $\theta^i=1$
and we use the saturation bound $\bar{c}^i=15$ for all
$i\in\until{6}$. Figure~\ref{fig::Ex3sim_Sat} shows the results of the
simulation for these two algorithms. Using high values for $\beta$ we
can reduce the tracking error, however, this results in larger driving
commands. As a result, both algorithms violate the saturation
bound. However, because the requirements of Lemma~\ref{lem::sat} are
satisfied in this example, as shown in Fig.~\ref{fig::Ex3sim_Sat}(b),
the ultimate tracking behavior of the agreement states of the
algorithm~\eqref{eq::DC2} despite the saturation resembles the
response of the algorithm~\eqref{eq::DC1} in the absence of saturation
bounds. There is not such guarantees for the algorithm~\eqref{eq::DC1}
(see Fig.~\ref{fig::Ex3sim_Sat}(a)).

\vspace{-6pt}
\section{Conclusions}\label{sec:conclusions}\vspace{-2pt}

This paper has addressed the multi-agent dynamic average consensus
problem over strongly connected and weight-balanced digraphs.  We have
proposed a distributed algorithm that makes individual agents track
the average of the dynamic inputs across the network with a
steady-state error. We have characterized how this error and the rate
of convergence depend on the design parameters of the proposed
algorithm, and identified special cases of inputs for which the
steady-state error is zero. Our algorithm enjoys the same convergence
properties in scenarios with time-varying topologies and is amenable
to discrete-time implementations.  We have also considered extensions
of the algorithm design that can handle limited control authority and
privacy preservation requirements against internal and external
adversaries.  Numerous avenues of research appear open for future
work, including the study of discrete-time implementations with the
features considered here (time-varying topologies, limited control
authority, and with privacy preservation features), the design of
provably-correct algorithms that do not require a priori
weight-balanced interaction topologies, and the application to
distributed estimation and map-merging scenarios.

\vspace*{-5ex}
\appendix
\section{Proof of the results of Section 4.3}\label{app:A}

Here, we provide the proof of the results presented in
Section~\ref{sec::Discrete}.

\begin{pro}[Proof of Lemma~\ref{lem::disc}]
  We can represent the zero-system of the discrete-time
  algorithm~\eqref{eq::DCDisc_smlr} in the following compact form
  \begin{align}\label{eq::UF1CDisc}
    \begin{bmatrix}
      \vect{x}(k+1)\\
      \vect{v}(k+1)
    \end{bmatrix}=\vect{P}_\delta\begin{bmatrix}
      \vect{x}(k)\\
      \vect{v}(k)
    \end{bmatrix},\quad \vect{P}_\delta=\vect{I}_{2N}+\delta\vect{A}.
  \end{align}
  where $\vect{A}$ is given in~\eqref{eq::UF1C}. Then,
  \begin{align*}
    \begin{bmatrix}
      \vect{x}(k)\\
      \vect{v}(k)
    \end{bmatrix}=\vect{P}_\delta^k\begin{bmatrix}
      \vect{x}(0)\\
      \vect{v}(0)
    \end{bmatrix}.
  \end{align*}
  In the proof of Lemma~\ref{lem::UF1} we showed that the eigenvalues
  of $\vect{A}$ are $-\alpha$ with multiplicity of $N$ and
  $-\beta\lambda_i$ for $i \in \until{N}$. Then, the eigenvalues of
  $\vect{P}_\delta$ are $1-\delta\alpha$ with multiplicity of $N$ and
  $1-\delta\beta\lambda_i$, where $i\in \until{N}$. Note that the
  eigenvalues of $\vect{I}_N-\delta\beta\vect{\mathsf{L}}$ are
  $1-\delta\beta\lambda_i$. Invoking~\cite[Lemma 3]{ROS-JAF-RMM:07},
  for a strongly connected and weight-balanced digraph, when
  $\delta\in(0,\min\{\alpha^{-1},\beta^{-1}(\text{d}_{\max}^{\text{out}})^{-1}\})$,
  the eigenvalues $1-\delta\beta\lambda_i$, $i=2,\dots,N$, are
  strictly inside the unit circle in the complex plane. Note that for
  $i=1$, $1-\delta\beta\lambda_i=1$. Therefore, we conclude that
    when
  $\delta\in(0,\min\{\alpha^{-1},\beta^{-1}(\text{d}_{\max}^{\text{out}})^{-1}\})$,
  for a strongly connected and weight-balanced digraph
  $\vect{P}_\delta$ has an eigenvalue equal to $1$ and the rest of the
  eigenvalues are located inside the unit circle. Therefore,
  $\vect{P}_\delta$ is a semi-convergent matrix, i.e.,
  $\lim_{k\to\infty}\vect{P}_\delta^k$ exists. Therefore
  \begin{displaymath}\begin{bmatrix}\vect{x}(k+1)\\
      \vect{v}(k+1)
    \end{bmatrix}-
    \begin{bmatrix}\vect{x}(k)\\
      \vect{v}(k)
    \end{bmatrix}\to\vect{0},~~\mathrm{as}~k\to\infty.\end{displaymath}
  Then,
  \begin{displaymath}\begin{bmatrix}\vect{x}(k+1)\\
      \vect{v}(k+1)
    \end{bmatrix}-
    \begin{bmatrix}\vect{x}(k)\\
      \vect{v}(k)
    \end{bmatrix}= \vect{P}_\delta\begin{bmatrix}\vect{x}(k)\\
      \vect{v}(k)\end{bmatrix}- \begin{bmatrix}\vect{x}(k)\\
      \vect{v}(k)\end{bmatrix}=
 \delta\vect{A}\begin{bmatrix}\vect{x}(k)\\
      \vect{v}(k)
    \end{bmatrix}\to\vect{0},~~\mathrm{as}~k\to\infty.\end{displaymath}
  As a result,
  \begin{equation}\label{eq::disc_limit}
    \lim_{k\to\infty}\left(\begin{bmatrix}\vect{x}(k)\\
        \vect{v}(k)
      \end{bmatrix}\right)=\mu\begin{bmatrix}\vect{1}_N\\
      -\alpha\vect{1}_N\end{bmatrix},\quad\mu\in\real.
  \end{equation}
  For a weight-balanced digraph, left multiplying the state equation
  of $\vect{v}$ by $\vect{1}^\top$, we obtain
  $\sum_{i=1}^Nv^{i}(k+1)=\sum_{i=1}^Nv^{i}(k)$. Consequently,
  $\sum_{i=1}^Nv^{i}(k)=\sum_{i=1}^Nv^{i}(0)$, $\forall
  \,k$. Invoking~\eqref{eq::disc_limit}, then at $k=\infty$ we have
  $-N\mu\alpha=\sum_{i=1}^Nv^{i}(0)$. As a result,
  $\mu=-\frac{\alpha^{-1}}{N}\SUM{i=1}{N}v^i(0)$.
\end{pro}
  
\begin{pro}[Proof of Theorem~\ref{thm::Disc_bound}]
  Consider the change of variables introduced in~\eqref{eq::xTOy},
  \eqref{eq::wTOv_disc}
  and~\eqref{eq::ChVarSim}. Then~\eqref{eq::DCDisc_smlr}, the
  equivalent representation of~\eqref{eq::DCDisc}, can be expressed in
  the following equivalent form
    \begin{align*}
      &
      \begin{bmatrix}
        p_1(k+1)\\
        q_1(k+1)
      \end{bmatrix}=\Tvect{P}_\delta\begin{bmatrix}
        p_1(k)\\
        q_1(k)
      \end{bmatrix},
          \\
      &
      \begin{bmatrix}\vect{p}_{2:N}(k+1)\\
        \vect{q}_{2:N}(k+1)\end{bmatrix}
      =\Bvect{P}_\delta\begin{bmatrix}\vect{p}_{2:N}(k)\\
        \vect{q}_{2:N}(k)\end{bmatrix}
      -
      \begin{bmatrix}
        \vect{0}
        \\
        \vect{R}^\top
      \end{bmatrix}
      (\Delta\vect{u}(k+1)-\Delta\vect{u}(k)+\delta\alpha
      \Delta\vect{u}(k)),
    \end{align*}
    where $\Tvect{P}_\delta=\vect{I}_2+\delta\Tvect{A}$ and
    $\Bvect{P}_\delta=\vect{I}_{N-2}+\delta\Bvect{A}$, with
    $\Tvect{A}$ and $\Bvect{A}$ are defined
    in~\eqref{eq::UF1C_Smlr}. For any given initial conditions, the
    solution of this difference equation is
    \begin{align*}
      p_1(k) =& p_1(0)-\delta
      \SUM{j=0}{k-1}(1-\delta\alpha)^jq_1(0), 
      \\
      q_1(k)=&(1-\delta\alpha)^k q_1(0), 
      \\
      \vect{p}_{2:N}(k) =& \vect{\Phi}(k,0) \vect{p}_{2:N}(0)-
      \SUM{j=0}{k-1}\vect{\Phi}(k-1,j)(1-\delta\alpha)^j(\vect{q}_{2:N}(0)
      + \vect{R}^\top\Delta\vect{u}(0))+ \nonumber
      \\
      &\SUM{j=0}{k-1}\vect{\Phi}(k-1,j)
      \vect{R}^\top\Delta\vect{u}(j),
      \\
      \vect{q}_{2:N}(k) =& (1-\delta\alpha)^k(\vect{q}_{2:N}(0) +
      \vect{R}^\top\Delta\vect{u}(0)) -
      \vect{R}^\top\Delta\vect{u}(k). 
    \end{align*}
  Recalling the change of variables~\eqref{eq::ChVarSim}, we have
  \begin{align*}
    \vect{y}(k)=&\vect{\mathsf{D}}_{11}\vect{y}(0) +
    \vect{\mathsf{D}}_{12}\vect{w}(0)-\vect{R}\SUM{j=0}{k-1}
    \vect{\Phi}(k-1,j)(1-\delta\alpha)^j\vect{R}^\top\Delta\vect{u}(0)+
    \vect{R}\SUM{j=0}{k-1}\vect{\Phi}(k-1,j)
    \vect{R}^\top\Delta\vect{u}(j),
  \end{align*}
  where
  \begin{align*}
    &\vect{\mathsf{D}}_{11} =
    (1-\alpha\delta\SUM{j=0}{k-1}(1-\delta\alpha)^j)\vect{r}
    \vect{r}^\top+\vect{R}\vect{\Phi}(k,0)\vect{R}^\top-\alpha
    \vect{R}(\SUM{j=0}{k-1}\vect{\Phi}(k-1,j)(1-\delta\alpha)^j)
    \vect{R}^\top,\nonumber
    \\
    &\vect{\mathsf{D}}_{12}=-\delta\SUM{j=0}{k-1}(1-\delta\alpha)^j\vect{r}\vect{r}^\top-
    \vect{R}(\SUM{j=0}{k-1}\vect{\Phi}(k-1,j)(1-\delta\alpha)^j)\vect{R}^\top.
  \end{align*}
  Because ${\vect{1}_N}^\top\vect{\Pi}_N=\vect{0}$,
  from~\eqref{eq::wTOv_disc} we can deduce
  $\SUM{i=1}{N}w^i(0)=\SUM{i=1}{N}v^i(0)$. Then, it is straightforward
  to obtain~\eqref{eq::xDCDisc_bound}.
\end{pro}

\begin{pro}[Proof of Corollary~\ref{cor::DCDiscProb1Sol}]
  We showed in Theorem~\ref{thm::Disc_bound} that, for any given
  stepsize, the bound~\eqref{eq::xDCDisc_bound} on the
  output $x^i$ of algorithm~\eqref{eq::DCDisc} holds. In the following, for
  the stepsizes satisfying $\delta \in(0,\min\{\alpha^{-1},
  \beta^{-1}(\text{d}_{\max}^{\text{out}})^{-1}\})$, we find the
  limiting value of the terms of this bound when $k\to\infty$.  Notice
  that $0<\delta<\alpha^{-1}$, then $0<(1-\alpha\delta)<1$. As a
  result, when $k\to\infty$ we have
  $\SUM{j=0}{k-1}(1-\delta\alpha)^j=(\delta\alpha)^{-1}$, leading to
  $(1-\alpha\delta\SUM{j=0}{k-1}(1-\delta\alpha)^j)\to0$ as
  $k\to\infty$. Recall $\vect{\Phi}(k,j) =
  (\vect{I}_{N-1}-\delta\beta\vect{R}^\top\vect{\mathsf{L}}\vect{R})^{k-j}$. Because
  $0<\delta<\beta^{-1}(\text{d}_{\max}^{\text{out}})^{-1}$, the
  spectral radius of $\vect{\Phi}(1,0)$ is less than one, therefore
  $\vect{\Phi}(k,0)\to0$ and $\SUM{j=0}{k-1}\vect{\Phi}(k-1,j)
  =(\delta\beta\vect{R}^\top\vect{\mathsf{L}}\vect{R})^{-1}$ as $k\to\infty$
  (see~\cite[Fact 10.3.1.xiii]{DSB:09}).  Also, there exists
  $\omega\in(0,1)$ such that $\rho(\vect{\Phi}(1,0))<\omega<1$. Then
  $\exists \,\mu>0$ such that $\Lnorm\vect{\Phi}(k-1,j)\Rnorm\leq
  \mu\, \omega^{k-1-j}$ for $0<j\leq k-1$,~\cite[pp. 26]{JPL:86}.  As
  a result, we have
  \begin{displaymath}
    \Lnorm\SUM{j=0}{k-1}\vect{\Phi}(k-1,j)(1-\delta\alpha)^j\Rnorm\leq
    \mu \SUM{j=0}{k-1}\omega^{k-1-j}(1-\delta\alpha)^j.
  \end{displaymath}
  Notice that
  \begin{displaymath}
    \SUM{j=0}{k-1}\omega^{k-1-j}(1-\delta\alpha)^j =
    \omega^{k-1}\SUM{j=0}{k-1}(\frac{1-\delta\alpha}{\omega})^j 
    = (1-\delta\alpha)^{k-1}\SUM{j=0}{k-1}(\frac{\omega}{1-\delta\alpha})^j.
  \end{displaymath}
  Then, as $k\to\infty$ we have
  \begin{displaymath}
    \Lnorm\SUM{j=0}{k-1}\!\!\vect{\Phi}(k-1,j)(1\!-\!\delta\alpha)^j
    \Rnorm\!\leq\!  
    \mu \!\SUM{j=0}{k-1}\!\omega^{k-1-j}(1\!-\!\delta\alpha)^j
    \!=
    \!\begin{cases}
      \mu\omega^{k-1}(1\!-\!\frac{1\!-\!\delta\alpha}{\omega})\!\!\to\!0,  
      &\omega\!>\!1\!-\!\delta\alpha,
      \\
     \mu (k\!-\!1)\omega^{k-1}\!\!\to\!0,&\omega\!=\!1\!-\!\delta\alpha,
      \\
      \mu(1\!-\!\delta\alpha)^{k-1} (1-\frac{\omega}{1\!-\!\delta\alpha})\!\!\to\!0,
      &\omega\!<\!1\!-\!\delta\alpha. 
    \end{cases}
  \end{displaymath}
  Invoking~\cite[Fact 8.18.12]{DSB:09}, we have
  \begin{displaymath}
    \Lnorm (\vect{R}^\top\vect{\mathsf{L}}\vect{R})^{-1}\Rnorm =
    \sigma_{\max}((\vect{R}^\top\vect{\mathsf{L}}\vect{R})^{-1})\leq
    \sigma_{\max}((\vect{R}^\top\Sym{\vect{\mathsf{L}}}\vect{R})^{-1}) = \Hlambda_2^{-1}.
  \end{displaymath}
  Also, notice that $\forall k\geq 0$, we have $\Lnorm \vect{R}^\top
  \Delta \vect{u}(k)\Rnorm =\Lnorm \vect{R}^\top \vect{\Pi}_N\Delta
  \vect{u}(k)\Rnorm\leq \Lnorm \vect{R}^\top\Rnorm\Lnorm
  \vect{\Pi}_N\Delta \vect{u}(k)\Rnorm \leq \gamma$.  Using the
  limiting values above, we can conclude that
  \begin{displaymath}
    \Lnorm\SUM{j=0}{k-1}\vect{\Phi}(k-1,j)
    \vect{R}^\top\Delta\vect{u}(j)\Rnorm\leq 
    \gamma \Lnorm\SUM{j=0}{k-1}\vect{\Phi}(k-1,j)\Rnorm = \gamma
    \Lnorm(\delta\beta\vect{R}^\top\vect{\mathsf{L}}\vect{R})^{-1}\Rnorm\leq
    \gamma/(\delta\beta\Hlambda_2). 
  \end{displaymath}
  This completes the proof. 
\end{pro}

\begin{pro}[Proof of Lemma~\ref{eq::Cnvrg_disc}]
  Using the change of variable~\eqref{eq::xTOy}, the algorithm~\eqref{eq::DCDisc} can be stated as follows
  (compact form)
  \begin{align}\label{eq::DCDisc_smlr2}
    \begin{bmatrix}
      \vect{y}(k+1)\\
      \vect{v}(k+1)
    \end{bmatrix}=\vect{P}_\delta \begin{bmatrix}
      \vect{y}(k)\\
      \vect{v}(k)
    \end{bmatrix}+\begin{bmatrix}
      \vect{\Pi}_N(\Delta\vect{u}(k)+\alpha\delta \vect{u}(k))\\
      \vect{0}
    \end{bmatrix},
  \end{align}
  where $\vect{P}_\delta$ is defined in~\eqref{eq::UF1CDisc}.  When
  condition (a) holds we have $\vect{\Pi}_N(\Delta
  \vect{u}(k)+\delta\alpha\vect{u}(k))\to\vect{0}$, as $k\to\infty$.
  Then~\eqref{eq::DCDisc_smlr2} is a linear system with a vanishing
  input $\vect{\Pi}_N(\Delta
  \vect{u}(k)+\delta\alpha\vect{u}(k))$. Therefore, it converges to
  the equilibrium of its zero-system. Notice that the system matrices
  of~\eqref{eq::DCDisc_smlr2} and~\eqref{eq::DCDisc_smlr} are the
  same. Therefore, when
  $\delta\in(0,\min\{\alpha^{-1},\beta^{-1}(\text{d}_{\max}^{\text{out}})^{-1}\})$,
  we can use result of Lemma~\ref{lem::disc} to conclude that
  $y^i(k)\to -\frac{\alpha^{-1}}{N}\SUM{i=1}{N}v^i(0)$, for $i \in
  \until{N}$.  Because $\SUM{i=1}{N}v^i(0)=0$, then we have $x^i(t)\to
  \avrg{u^j(k)}$ globally asymptotically for $i \in \until{N}$.  Next,
  notice that using the change of variables~\eqref{eq::xTOy}
  and~\eqref{eq::wTOv_disc} another equivalent representation
  of~\eqref{eq::DCDisc} can be stated as follows 
  \begin{align}\label{eq::DCDisc_smlr3}
    \begin{bmatrix}
      \vect{y}(k+1)\\
      \vect{w}(k+1)
    \end{bmatrix}=\vect{P}_\delta \begin{bmatrix}
      \vect{y}(k)\\
      \vect{w}(k)
    \end{bmatrix}-\begin{bmatrix}
      \vect{0}\\
      \vect{\Pi}_N(\Delta\vect{u}(k+1)-\Delta\vect{u}(k)+\delta\alpha
      \Delta\vect{u}(k))
      \\
    \end{bmatrix},
  \end{align}
  where $\vect{P}_\delta$ again is defined in~\eqref{eq::UF1CDisc}.
  When condition (b) holds we have $\vect{\Pi}_N(\Delta
  \vect{u}(k)+\delta\alpha\vect{u}(k))\to\vect{0}$ as $k\to\infty$.
  Then,~\eqref{eq::DCDisc_smlr3} is a linear system with a vanishing
  input
  $\vect{\Pi}_N(\Delta\vect{u}(k+1)-\Delta\vect{u}(k)+\delta\alpha
  \Delta\vect{u}(k))$. Then, using a similar argument used
  for~\eqref{eq::DCDisc_smlr2} above, we can show that
  in~\eqref{eq::DCDisc_smlr3} $y^i(k)\to
  -\frac{\alpha^{-1}}{N}\SUM{i=1}{N}w^i(0)$ as $k\to\infty$ for $i \in
  \until{N}$. Using~\eqref{eq::wTOv_disc}, we can show
  $\SUM{i=1}{N}w^i(0)=\SUM{i=1}{N}v^i(0)$. As a result $x^i(k)\to
  \avrg{u^j(k)}$ globally asymptotically for $i \in \until{N}$.
\end{pro}

\section{Supporting material for the proof of
  Lemma~\ref{lem::sat}}\label{app:B}

The following results are used in the proof of Lemma~\ref{lem::sat}.

\begin{pro}\label{pro::sat1}
  Consider the following system where $x,w,\beta\in\real$, $\beta>0$
  and $w$ is a piece-wise continuous time-varying signal
  \begin{equation}\label{eq::sys_sat_eval1}
    \dot{y}=-\beta\sat{\bar{c}}{y-w}-\beta w.
  \end{equation}
 Assume that $||w||_{\text{ess}}<
  \bar{c}$. Then, for any initial condition $y(0)\in\real$, $y(t)\to
  0$ asymptotically.
\end{pro}
\begin{pro}
  Consider the candidate Lyapunov function $V=\frac{1}{2\beta}y^2$
  with derivative $\dot{V}=-y\sat{\bar{c}}{y-w}-yw$ along the
  trajectories of~\eqref{eq::sys_sat_eval1}.  To prove that $\dot{V}$
  is negative definite, first note that because $||w||_{\text{ess}}<
  \bar{c}$, we have that if $y-w>\bar{c}$ then $y>\bar{c}+w>0$ and if
  $y-w<-\bar{c}$ then $y<-\bar{c}+w<0$.  As a result,
  \begin{align*}
    \dot{V} = 
    \begin{cases}
      -y(\bar{c}+w)\leq -(\bar{c}-||w||_{\text{ess}})|y|<0, &
      ~\text{if} ~y-w>\bar{c},
      \\
      -y^2<0, & ~\text{if} ~|y-w|\leq\bar{c},
      \\
      -y(-\bar{c}+w)\leq -(\bar{c}-||w||_{\text{ess}})|y|<0 ,
      &~\text{if}~ y-w<-\bar{c} .
    \end{cases}
  \end{align*}
  All the conditions of the Lyapunov stability analysis for
  non-autonomous systems \cite[Theorem 4.9]{HKK:02} are satisfied
  globally. Therefore, $y(t)\to0$ globally asymptotically as
  $t\to\infty$.
\end{pro}

\begin{pro}\label{pro::sat2}
  Consider the following system where $x,u\in\real$ and $u$ is a piece-wise continuous time-varying
  signal, 
  \begin{equation}\label{eq::sys_sat_eval2}
    \dot{x}=-\sat{\bar{c}}{\beta(x-u)-\dot{u}} ,
  \end{equation}
  Assume $u$ and its derivative $\dot{u}$ are both
  essentially bounded signals, and there is some finite $t^\star>0$
  such that for all $t\geq t^\star$, $|\dot{u}(t)|< \bar{c}$. Then,
  for any initial condition $x(0)\in \real$ we have $x(t)\to u(t)$
  asymptotically.
\end{pro}
\begin{pro}
  Given that~\eqref{eq::sys_sat_eval2} is ISS, c.f.~\cite{EDS:95a},
  and since $\beta u+\dot{u}$ is bounded, for any finite initial
  condition $x(0)$, there is a finite $\mu(x(0))>0$ such that we have
  $|x|<\mu(x(0))$ for all $t\geq0$.  Under the change of variables $
  y=\beta(x-u)$, equation~\eqref{eq::sys_sat_eval2} can be written in
  the following equivalent form
  \begin{equation}\label{eq::sys_sat_eval3}
    \dot{y}=-\beta\sat{\bar{c}}{y-\dot{u}}-\beta \dot{u} .
  \end{equation}
  Since the solutions of~\eqref{eq::sys_sat_eval2} are all bounded and
  because both $u$ and $x$ are bounded signals, starting from any
  initial condition, we have the guarantee that the solutions
  of~\eqref{eq::sys_sat_eval3}) are also bounded. Since the input
  $\dot{u}$ to the system~\eqref{eq::sys_sat_eval3} satisfies the
  conditions of Proposition~\ref{pro::sat1} after some finite time
  $t^\star$, we can conclude that $y(t)\to0$, or equivalently $x(t)\to
  u(t)$, globally asymptotically.
\end{pro}


\vspace*{-3ex}
\begin{thebibliography}{10}
\providecommand{\url}[1]{\texttt{#1}}
\providecommand{\urlprefix}{URL }
\expandafter\ifx\csname urlstyle\endcsname\relax
  \providecommand{\doi}[1]{doi:\discretionary{}{}{}#1}\else
  \providecommand{\doi}{doi:\discretionary{}{}{}\begingroup
  \urlstyle{rm}\Url}\fi

\bibitem{PY-RAF-KML:08}
Yang P, Freeman R, Lynch K. Multi-agent coordination by decentralized
  estimation and control. \emph{IEEE Transactions on Automatic Control}  2008;
  \textbf{53}(11):2480Ð--2496.

\bibitem{SM:07e}
Mart{\'\i}nez S. Distributed interpolation schemes for field estimation by
  mobile sensor networks. \emph{IEEE Transactions on Control Systems
  Technology}  2010; \textbf{18}(2):491--500.

\bibitem{JC:07-tac}
Cort\'es J. Distributed {K}riged {K}alman filter for spatial estimation.
  \emph{IEEE Transactions on Automatic Control}  2009;
  \textbf{54}(12):2816--2827.

\bibitem{ROS-JSS:05}
Olfati-Saber R, Shamma J. Consensus filters for sensor networks and distributed
  sensor fusion. \emph{{IEEE} Int. Conf. on Decision and Control and European
  Control Conference}, Seville, Spain, 2005; 6698--6703.

\bibitem{ROS:07}
Olfati-Saber R. Distributed {K}alman filtering for sensor networks.
  \emph{{IEEE} Int. Conf. on Decision and Control}, New Orleans, LA, 2007;
  5492--5498.

\bibitem{RA-JC-CS:10j}
Arag\"u\'es R, Cort\'es J, Sag\"u\'es C. Distributed consensus on robot
  networks for dynamically merging feature-based maps. \emph{IEEE Transactions
  on Robotics}  2012; \textbf{28}(4):840--854.

\bibitem{PY-RAF-KML:07}
Yang P, Freeman R, Lynch K. Distributed cooperative active sensing using
  consensus filters. \emph{{IEEE} Int. Conf. on Robotics and Automation}, Roma,
  Italy, 2007; 405--410.

\bibitem{ROS-JAF-RMM:07}
Olfati-Saber R, Fax JA, Murray RM. Consensus and cooperation in networked
  multi-agent systems. \emph{Proceedings of the IEEE}  2007;
  \textbf{95}(1):215--233.

\bibitem{WR-RWB:08}
Ren W, Beard RW. \emph{Distributed Consensus in Multi-vehicle Cooperative
  Control}. Communications and Control Engineering, Springer, 2008.

\bibitem{WR-YC:11}
Ren W, Cao Y. \emph{Distributed Coordination of Multi-Agent Networks}.
  Communications and Control Engineering, Springer: New York, 2011.

\bibitem{MM-ME:10}
Mesbahi M, Egerstedt M. \emph{Graph Theoretic Methods in Multiagent Networks}.
  Princeton Series in Applied Mathematics, Princeton University Press, 2010.

\bibitem{FB-JC-SM:09}
Bullo F, Cort{\'e}s J, Mart{\'\i}nez S. \emph{Distributed Control of Robotic
  Networks}. American Mathematical Society, Princeton University Press, 2009.
  Available at http://www.coordinationbook.info.

\bibitem{DPS-ROS-RMM:05b}
Spanos D, Olfati-Saber R, Murray R. Dynamic consensus on mobile networks.
  \emph{{IFAC} {W}orld {C}ongress}, Prague, Czech Republic, 2005.

\bibitem{ROS-RMM:04}
Olfati-Saber R, Murray R. Consensus problems in networks of agents with
  switching topology and time-delays. \emph{IEEE Transactions on Automatic
  Control}  September 2004; :1520--1533.

\bibitem{RAF-PY-KML:06}
Freeman RA, Yang P, Lynch KM. Stability and convergence properties of dynamic
  average consensus estimators. \emph{{IEEE} Int. Conf. on Decision and
  Control}, San Diego, CA, 2006; 398--403.

\bibitem{HB-RAF-KML:10}
Bai H, Freeman R, Lynch K. Robust dynamic average consensus of time-varying
  inputs. \emph{{IEEE} Int. Conf. on Decision and Control}, Atlanta, GA, USA,
  2010; 3104--3109.

\bibitem{MZ-SM:08a}
Zhu M, Mart{\'\i}nez S. Discrete-time dynamic average consensus.
  \emph{Automatica}  2010; \textbf{46}(2):322--329.

\bibitem{WR:07a}
Ren W. Multi-vehicle consensus with a time-varying reference state.
  \emph{Systems and Control Letters}  2007; \textbf{56}(2):474--483.

\bibitem{GS-YH-KHJ:12}
Shi G, Hong Y, Johansson KH. Connectivity and set tracking of multi-agent
  systems guided by multiple moving leaders. \emph{IEEE Transactions on
  Automatic Control}  2012; \textbf{57}(3):663--676.

\bibitem{GS-KHJ:13}
Shi G, Johansson KH. Robust consensus for continuous-time multiagent dynamics.
  \emph{SIAM Journal on Control and Optimization}  2013;
  \textbf{51}(5):3673--3691.

\bibitem{JH:04}
Hespanha J. Uniform stability of switched linear systems: extensions of
  lasalle's invariance principle. \emph{IEEE Transactions on Automatic Control}
   2004; \textbf{49}(4):470--482.

\bibitem{DSB:09}
Bernstein D. \emph{Matrix {M}athematics: theory, facts, and formulas}. 2 edn.,
  Princeton University Press, 2009.

\bibitem{EDS-YW:95}
Sontag ED, Wang Y. On characterizations of the input-to-state stability
  property. \emph{Systems and Control Letters}  1995; \textbf{24}(5):351--359.

\bibitem{LV-DC-DL:05}
Vu L, Chatterjee D, Liberzon D. I{S}{S} of switched systems and applications to
  switching adaptive control. \emph{{IEEE} Int. Conf. on Decision and Control},
  Seville, Spain, 2005; 120--125.

\bibitem{JPL:86}
LaSalle JP. \emph{The Stability and Control of Discrete Processes},
  \emph{Applied Mathematical Sciences}, vol.~62. Springer, 1986.

\bibitem{HKK:02}
Khalil HK. \emph{Nonlinear Systems}. 3 edn., Prentice Hall, 2002.

\bibitem{EDS:95a}
Sontag ED. On the input-to-state stability property. \emph{European Journal of
  Control}  1995; \textbf{1}:24--36.

\end{thebibliography}
\end{document}